\input ./style/arxiv-general.cfg
\documentclass[aap,MSNbibl,seceqn,dvips]{arximspdf}
\makeatletter
   \@ifpackageloaded{graphicx}{}{\usepackage{graphicx}}
\makeatother
\usepackage{mathbh}
\usepackage{url,breakurl}

\doi{10.1214/14-AAP1053} 
\volume{25}
\issue{5}
\pubyear{2015}
\firstpage{2462}
\lastpage{2502}
\docsubty{FLA}

\makeatletter
\newcommand{\rrvert}{\vert}
\newcommand{\llvert}{\vert}
\newcommand{\eqref}[1]{(\ref{#1})}

\newtheorem{Theorem}{Theorem}[section]
\newtheorem{Lemma}[Theorem]{Lemma}
\newtheorem{Proposition}[Theorem]{Proposition}

\newtheorem{Corollary}[Theorem]{Corollary}

\newproclaim{Remark}[Theorem]{Remark}
\newproclaim{Example}[Theorem]{Example}
\newproclaim{OP}[Theorem]{Open Problem}

\newproclaim{Definition}[Theorem]{Definition}
\newproclaim{Construction}[Theorem]{Construction}

\newcommand{\cH}{\mathcal{H}}
\newcommand{\cS}{\mathcal{S}}

\newcommand{\HH}{\mathcal{H}}
\newcommand{\prob}{\mathbb{P}}
\newcommand{\GG}{\mathcal{G}}
\newcommand{\BP}{\mathsf{BP}}
\newcommand{\Yu}{\mathsf{Yu}}
\newcommand{\card}{\operatorname{card}}
\newcommand{\expec}{\mathbb{E}}
\def\ind{\mathbh{1}}
\newcommand{\sss}{}
\newcommand{\convp}{\stackrel{\sss{\mathbb P}}{\longrightarrow}}

\newcommand{\mathS}{\mathcal{S}} 

\newcommand{\pa}{preferential attachment }

\newcommand{\lamb}{\operatorname{Lam}}

\makeatother

\begin{document}
\begin{frontmatter}

\title{Twitter event networks and the Superstar model}
\runtitle{Twitter event networks and the Superstar model}

\begin{aug}
\author[A]{\fnms{Shankar}~\snm{Bhamidi}\corref{}\thanksref{T1}\ead[label=e1]{bhamidi@email.unc.edu}},
\author[B]{\fnms{J.~Michael}~\snm{Steele}\ead[label=e2]{steele@wharton.upenn.edu}}
\and
\author[C]{\fnms{Tauhid}~\snm{Zaman}\ead[label=e3]{zlisto@mit.edu}}
\runauthor{S. Bhamidi, J.~M. Steele and T. Zaman}
\thankstext{T1}{Supported in part by NSF-DMS Grants 1105581 and 1310002
and SES-1357622.}
\affiliation{University of North Carolina, University of Pennsylvania and Massachusetts Institute of Technology}
\address[A]{S. Bhamidi\\
Department of Statistics\\
\quad and Operations Research\\
University of North Carolina\\
Chapel Hill, North Carolina 27599\\
USA\\
\printead{e1}} 
\address[B]{J. M. Steele\\
The Wharton School\\
Department of Statistics\\
Huntsman Hall 447\\
University of Pennsylvania\\
Philadelphia, Pennsylvania 19104\\
USA\\
\printead{e2}}
\address[C]{T. Zaman\\
Sloan School of Management\\
Massachusetts Institute of Technology\\
Cambridge, Massachusetts 02142\\
USA\\
\printead{e3}}
\end{aug}

\received{\smonth{11} \syear{2012}}
\revised{\smonth{7} \syear{2014}}

%
\begin{abstract}
Condensation phenomenon is often observed
in social networks such as Twitter where
one ``superstar'' vertex gains a positive fraction of the
edges, while the remaining empirical degree distribution
still exhibits a power law tail. We formulate a mathematically
tractable model for this phenomenon that provides a better fit
to empirical data than the standard preferential attachment model
across an array of networks observed in Twitter. Using embeddings
in an equivalent continuous time version of the process, and
adapting techniques from the stable age-distribution theory of
branching processes, we prove limit results for the proportion
of edges that condense around the superstar, the degree distribution
of the remaining vertices, maximal nonsuperstar degree asymptotics
and height of these random trees in the large network limit.
\end{abstract}

%
\end{frontmatter}
%
\section{Retweet graphs and a mathematically tractable model}\label{sec:intro}

Our goal here is to provide a simple model that captures the most
salient features of a natural graph that is determined by the
Twitter traffic generated by public events. In the Twitter world (or
Twitterverse), each user has a set of followers; these are people who
have signed-up to receive the tweets of the user.
Here, our focus is on \emph{retweets}; these are tweets by a user who forwards
a tweet that was received from another user. A~retweet is sometimes
accompanied with comments by the retweeter.

Let us first start with an empirical example that contains all the
characteristics observed in a wide array of such retweet networks. Data
was collected during the Black Entertainment Television (BET) Awards of 2010.
We first considered all tweets in the Twitterverse that were posted
between 10 AM and 4 PM (GMT) on the day of the ceremony, and we then restricted
attention to all the tweets in the Twitterverse that contained the term
``BET Awards.'' We view the posters of
these tweets as the vertices of an undirected simple graph where there
is an edge between vertices
$v$ and $w$ if $w$ retweets a tweet received from $v$, or vice-versa.
We call this graph the \emph{retweet graph.}

%
\begin{figure}

\includegraphics{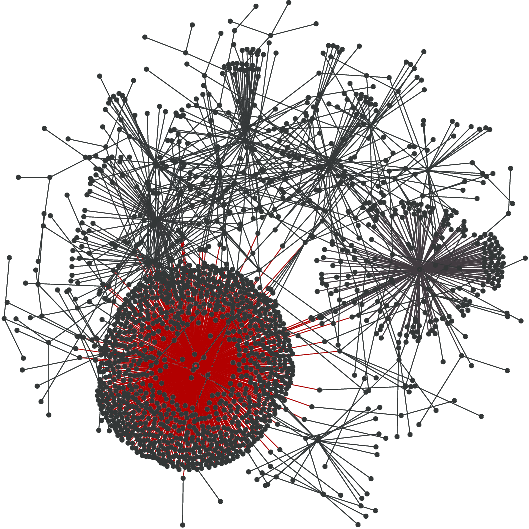}

\caption{Giant component of the 2010 BET Awards retweet graph.}
\label{fig:betawards}
\end{figure}

In the retweet graph for the 2010 BET Awards, one finds a single giant
component (see Figure~\ref{fig:betawards}).
There are also many small components
(with five or fewer vertices)\vadjust{\goodbreak}
and a large number of isolated vertices. The giant component is also
approximately a tree in the sense that if we remove 91 edges from the
graph of 1724 vertices and 1814 edges we obtain an honest tree.
Finally, the most compelling feature of this empirical tree is that it
has one vertex of \emph{exceptionally} large degree.
This ``superstar'' vertex has degree 992, so it is connected to more
than 57\% of the vertices.
As it happens, this ``superstar'' vertex corresponds to the
pop-celebrity Lady Gaga who received an award at the ceremony.

%
%
\subsection{Superstar model for the giant component}
Our main observation is that the qualitative and quantitative features
{{of}} the giant component {{in a wide array of}}
retweet graphs may be captured rather well by a simple one-parameter model.
The construction of the model only makes an obvious modification of the
now classic preferential attachment model,
but this modification turns out to have
richer consequences than its simplicity would suggest. Naturally, the
model has the ``superstar'' property baked into the cake,
but a surprising consequence is that the distribution of the degrees of
the nonsuperstar vertices is {{quite}} different
from what one finds in the preferential attachment model.

Our model is a graph evolution process that we denote by $\{G_n$,
$n=1,2, \ldots\}$. The graph $G_1$ consists of
the single vertex $v_0$, that we call the superstar. The graph $G_2$
then consists of the superstar $v_0$, a nonsuperstar $v_1$,
and an edge between the two vertices. For $n\geq2$, we construct
$G_{n+1}$ from $G_{n}$ by attaching the vertex $v_{n}$ to the superstar
{{$v_0$}} with probability $0< p< 1$
while with probability $q=1-p$ we attach $v_{n}$ to a nonsuperstar
according to the classical \emph{preferential attachment rule}.
That is, with probability $q$ the nonsuperstar
$v_{n}$ is attached to one of the nonsuperstars $\{v_1, v_2,\ldots,
v_{n-1}\}$ with probability that is proportional to the degree of $v_i$
in~$G_n$.

%
%
\subsection{Organization of the paper}
In the next section, we state the main results for the Superstar
model. In Section~\ref{sec:discussion}, we consider previous work on
Twitter networks and explore the connection between our model and
existing models. {{In this section, we also describe two variants of
the basic Superstar model (linear attachment and uniform attachment)
that can be rigorously analyzed using the same mathematical methodology
developed in this paper.}} In Section~\ref{sec:data}, we study the
performance of this model on various real networks constructed from the
Twitterverse and we compare our model to the standard preferential
attachment model. Section~\ref{sec:proofs} is the heart of the paper.
Here, we construct a special two-type continuous time branching process
that turns out to be equivalent to the Superstar model and analyze
various structural properties of this continuous time model. In
Section~\ref{sec:equiv}, we prove the equivalence between the
continuous time
model and the Superstar model through a \emph{surgery} operation. In
Section~\ref{sec:comp-proof}, we complete the proofs of all the main results.

\section{Mathematical results for the Superstar model}\label{sec:results}

Let $\{G_n$, $n=1,2, \ldots\}$ denote {{the}} graph process that
evolves according to the Superstar model with parameter $0 < p < 1$. We
shall think about all the processes constructed on a single probability
space through the obvious sequential growth mechanism so that one can
make almost sure statements.
The degree of the vertex
$v$ in the graph $G$ is denoted by $ \deg(v, G)$. The first result
describes asymptotics of the condensation phenomenon around the
superstar. The result is an immediate consequence of the {{definition
of the model and the strong law of large numbers. Since}} it
is a defining element of our model, we set the result out as a theorem.

%
\begin{Theorem}[(Superstar strong law)] \label{thm:SSstronglaw}
With probability one, we have
%
%
\begin{equation}
\label{SSLLN} \lim_{n \rightarrow\infty} \frac{1}{n}
\deg(v_0, G_n) = p.
\end{equation}
\end{Theorem}
The next result describes the asymptotic degree distribution.
%
%
\begin{Theorem}[(Degree distribution strong law)]\label{thm:DegreeDist}
With probability one, we have
\[
\lim_{n \rightarrow\infty} \frac{1}{n} \card\bigl\{1\leq j\leq
n\dvtx
\deg(v_j, G_n) = k\bigr\} = \nu_{\mathrm{SM}} (k,p ),
\]
where $\nu_{\mathrm{SM}}(\cdot, p)$ is the probability mass function defined on
$ \{1,2,\ldots\}$ by
\[
\nu_{\mathrm{SM}} (k,p )=\frac{2-p}{1-p} (k-1 )!\prod
_{i=1}^k \biggl(i+\frac{2-p}{1-p}
\biggr)^{-1}.
\]
\end{Theorem}
%
%
\begin{Remark}
One should note that the above theorem implies that the degree
distribution of the nonsuperstar vertices has a power law tail. Specifically,
\[
\frac{2-p}{1-p} (k-1 )!\prod_{i=1}^k
\biggl(i+\frac
{2-p}{1-p} \biggr)^{-1}\sim C_pk^{-\beta}\qquad
\mbox{as $k\to\infty$},
\]
for the constants
\[
\beta= 3+p/(1-p), \qquad{{C_p = \biggl(\frac{2-p}{1-p}
\biggr)^2\Gamma\biggl(\frac{2-p}{1-p} \biggr)}},
\]
where $\Gamma(x)$ is the gamma function. This should be contrasted with
the standard preferential attachment model {{(with no superstar
attachment)}} whose degree distribution scales like $k^{-3}$ as $k\to
\infty$. Thus, although one might expect that this variation in the
attachment scheme implies that a fraction $1-p$ of the vertices still
continue to {{perform}} preferential attachment, and thus the degree
distribution should still have a power law exponent of $3$; in reality,
this attachment scheme has a major effect on the degree distribution.
{{One}} requires a careful analysis of the different time-scales of
the associated continuous time branching process to tease out
asymptotic properties of the model.
\end{Remark}

The next theorem concerns the largest degree amongst all the
nonsuperstar vertices $ \{v_i\dvtx1\leq i\leq n \}$. Let
\[
\Upsilon_n:= \max_{1\leq i\leq n} \operatorname{deg}(v_i,
G_n).
\]
%
%
\begin{Theorem}[(Maximal nonsuperstar degree)]\label{thm:DegreeEvolutionLaw}
Let $\gamma= (1-p)/(2-p)$. There exists a random variable $\Delta^*$
with $\prob(0< \Delta^*< \infty) =1$ such that
\[
\lim_{n \rightarrow\infty}\frac{1}{n^\gamma}\Upsilon_n \stackrel
{\prob} {\longrightarrow} \Delta^*.
\]
%
\end{Theorem}

The almost sure linear growth of the degree of the superstar
(Theorem~\ref{thm:SSstronglaw}) is endemic to our construction. For standard
preferential attachment (with no superstar attachment mechanism), the
maximal degree grows like $\Theta_P(n^{1/2})$ (cf. \cite
{mori2007degree}). {{Thus, the superstar attachment affects the scaling
of the maximal degree as well. }}

Recall that $G_n$ is a tree. View this tree as rooted at the superstar
{{vertex}} $v_0$. {{Write}} $\HH(G_n)$ for the graph distance of
the vertex furthest from the root. {{Thus,}} $\cH(G_n)$ is the
height of {{the random tree}} $G_n$. Theorem~\ref{thm:SSstronglaw}
implies that a fraction $p$ of the vertices in the network are directly
connected to the superstar. One {{might}} wonder if this reflects a
general property of the network, {{namely does}} $\HH(G_n) = O_p(1)$
as $n\to\infty$? The next theorem shows that in fact the height of the
tree increases logarithmically in the size of the network. Let $\lamb
(\cdot)$ be the Lambert special function (cf. \cite
{corless1996lambertw}) and recall that $ \lamb(1/e)\approx0.2784$.
%
%
\begin{Theorem}[(Logarithmic height scaling)]\label{thm:height}
With probability one, we have
\[
\lim_{n \rightarrow\infty}\frac{1}{\log{n}}\HH(G_n) =
\frac
{1-p}{{{\lamb}}(1/e) (2-p)}.
\]
\end{Theorem}

\section{Related results and questions}
\label{sec:discussion}
In this section, we briefly discuss the connections between this model
and some of the more standard models in the literature as well as
extensions of the results in the paper. We also discuss previous
empirical research done on the
structure of Twitter networks.
%
\subsection{Preferential attachment}
This has become one of the standard\break workhorses in the complex networks
community. It is {{well nigh impossible to compile even a
representative list of references}}; see \cite{szymanski1987nonuniform}
where it was introduced in the combinatorics community, \cite
{barabasi1999emergence} for bringing this model to the attention of the
networks community, \cite{newman2003structure,dorogovtsev2002evolution} for survey level treatments of a wide array
of models, \cite{Bollobas:2001:DSS:379831.379835} for the first
rigorous results on the asymptotic degree distribution and \cite
{cooper2003general,bollobas2003mathematical,rudas-1}
and~\cite{durrett-rg-book} and the references therein for more general
models and results.
{{Let us briefly describe the simplest model in this class of models.}}
One starts with two vertices connected by a single edge as in
the Superstar model. {{Then}} each new vertex joins the system by
connecting to a single vertex in the current tree by choosing this
{{extant}} vertex with probability proportional to its {{current}}
degree. In this case, one can show \cite
{Bollobas:2001:DSS:379831.379835} that there exists a limiting
asymptotic degree distribution, {{namely}} with probability one
\[
\lim_{n \rightarrow\infty} \frac{1}{n} \operatorname{card} \bigl
\{ 1
\leq j\leq n\dvtx\deg(v_j, G_n)=k \bigr\}=
\frac{4}{k(k+1)(k+2)}.
\]
{{Thus, the asymptotic degree distribution exhibits}} a degree
exponent of three. The Superstar model changes the degree exponent of
the nonsuperstar vertices from three to $(3-2p)/(1-p)$ (see
Theorem~\ref
{thm:DegreeEvolutionLaw}).
Further, for the preferential attachment model, the maximal degree
scales like
$n^{1/2}$ \cite{mori2007degree}, while for the Superstar model, the
maximal nonsuperstar degree scales like $n^{\gamma}$ with $\gamma=
(1-p)/(2-p)$.
\subsection{Statistical estimation} We use real data on various
Twitter streams to analyze the empirical performance of the Superstar
model and compare this with typical preferential attachment models in
Section~\ref{sec:data}. Estimating the parameters from the data raises
a host of new interesting statistical questions. See \cite
{wiuf2006likelihood} where such questions were first raised and
likelihood based schemes were proposed in the context of usual
preferential attachment models. Considering how often such models are
used to draw quantitative conclusions about real networks, proving
consistency of such procedures as well as developing methodology to
compare different estimators in the context of models of evolving
networks would be of great interest to a number of different fields.
\subsection{Stable age distribution} The proofs for the degree
distribution build heavily on the analysis of the stable age
distribution for a single type continuous time branching process in
\cite{nerman1981convergence}. We extend this analysis to the context of
a two-type variant whose evolution mirrors the discrete type model.
Using Perron--Frobenius theory, a wide array of structural properties
are known about such models (see~\cite{jagers-nerman-comp-multi}). The
models used in our proof technique are relatively simpler and we can
give complete proofs using special properties of the continuous time
embeddings, including special martingales that play an integral role in
the treatment (see, e.g., Proposition~\ref{pro:conts-mp}). There have
been a number of recent studies on various preferential attachment
models using continuous time branching processes; see, for example,
\cite{rudas2007random,athreya2008growth,deijfen2010random}. For the
usual preferential attachment model ($p=0$), \cite{pittel-height}
{{uses}} embeddings in continuous time and results on the first birth
time in such branching processes (see \cite{kingman-age-depen}) to show
that the height satisfies
\[
\frac{\HH(\GG_n)}{\log{n}} \stackrel{\mathrm{a.s.}}{\longrightarrow}\frac
{1}{2\lamb(1/e)}.
\]
Here, we use a similar technique, but we first need to extend \cite
{kingman-age-depen} to the setting of multitype branching processes.
%
\subsection{Previous analysis of Twitter networks} The majority
of work analyzing Twitter networks has been empirical in nature.
In one of the earliest studies of Twitter networks \cite
{ref:twitterMoon}, the
authors looked at the degree distribution of the different networks in Twitter,
including retweet networks associated with individual topics. Power-laws
were observed, but no model was proposed to describe the network evolution.
In \cite{asur2011trends}, the link between maximum degree
and the range of time for which a topic was popular or ``trending''
was investigated. Correlations between the degree in
retweet graphs and the Twitter follower graph for different
users was studied in \cite{ref:miacha}. These empirical
analyses provided many important insights into the structure
of networks in Twitter. However, the lack of a model
to describe the evolution of these networks is one of the important unanswered
questions in this field, and the rigorous analysis of such a model
has not {{yet been considered}}. Our work here presents one of the first
such models that produces predictions that match Twitter data and also
{{provides a rigorous theoretical analysis of the proposed model}}.

\subsection{Related models} One of the main aims of this work is to
{{develop mathematical techniques that extend in a straightforward
fashion to variants of the Superstar model}}. We state results for two
such models in this section. {{We will describe how to extend the
proofs for the Superstar model to these variants in Section~\ref
{sec:extension-proofs}.}} We first start with {{the}} \emph{superstar
linear preferential attachment}.
%
Fix a parameter $a > -1$. The linear preferential attachment model is
described as follows: As before new vertices attach to vertex $v_0$
with probability $p$. With probability $q:=1-p$ the new vertex attaches
to one of the extant vertices $v$, with probability proportional to the
$d(v)+a$ where $d(v)$ is the present degree of the vertex. As before,
by construction the degree of the superstar $v_0$ scales like $\sim pn$
as $n\to\infty$. The techniques in the paper extend with simple
modifications to prove the following.

%
\begin{Theorem}[(Linear superstar preferential attachment)]
\label{thm:superstar-linear}
Fix $a> -1$ and $p\in(0,1)$. In the linear Superstar model, one has
for all $k\geq1$, with probability one
\begin{eqnarray*}
&&\lim_{n\to\infty} {{\frac{1}{n}{\card\bigl\{1\leq j\leq n
\dvtx\deg(v_j, G_n) = k \bigr\}}}}\\
&&\qquad = \frac{2-p+a}{1-p}
\frac{\prod_{j=1}^{k-1}(j+a)}{\prod_{i=1}^k(i+({(2-p)}/{(1-p)})(1+a) )}.
\end{eqnarray*}
Further, for $\gamma(a) = (1-p)/(2-p+a)$, there exists a random
variable $0< \Delta^*(a) < \infty$ a.s. such that the largest degree
other than the superstar satisfies
\[
n^{-\gamma(a)}\max_{ \{1\leq i\leq n \}} \deg(v_i) \convp
\Delta^*(a) \qquad\mbox{as } n\to\infty.
\]
\end{Theorem}
Similarly, one can show that the height of the linear Superstar model
scales like $\kappa(a)\log{n}$ for a limit constant $0 < \kappa(a) <
\infty$.

We next consider the case of the less realistic Superstar model with
\emph{uniform attachment}. Here, each new vertex attaches to the
superstar $v_0$ with probability $p$ or to one of the remaining
vertices uniformly at random (irrespective of the degree). Although
less realistic in the context of social networks, this is the superstar
variant of the random recursive tree {{a model of a growing tree where
each new vertex attaches to a uniformly chosen extant vertex}}.
The random recursive tree has been a model of great interest in the
combinatorics and computer science community (see the survey \cite
{smythe1995survey}). This model differs from the previous models with
the limiting degree distribution possessing exponential tails while the
maximal degree only growing logarithmically in the size of the network.

%
\begin{Theorem}[(Superstar uniform attachment)]
\label{thm:superstar-uniform}
{{Let $q:=1-p$.}} For the uniform attachment model, one has for all
$k\geq1$ that with probability one
\[
\lim_{n\to\infty} \frac{1}{n} \card{\bigl\{1\leq j\leq n\dvtx
\deg(v_j, G_n) = k\bigr\}} = \frac{1}{1+q} \biggl(
\frac{q}{1+q} \biggr)^{k-1},
\]
and the maximal nonsuperstar degree satisfies
\[
\lim_{n\to\infty}\frac{\max_{1\leq i\leq n} \deg(v_i)}{\log{n}}
\convp\frac{1}{\log{{(1+q)}/{q}}}.
\]
\end{Theorem}

\section{Retweet graphs for different public events}\label{sec:data}

We collected tweets from the Twitter firehose for thirteen different
public events, such as sports matches and musical performances \cite{ref:msr}.
The Twitter firehose
is the full feed of all public tweets that is accessed via Twitter's Streaming
Application Programming Interface~\cite{ref:firehose}. By using the
Twitter firehose, we were able to access all public tweets in the Twitterverse.

For each public event $E\in\{1,2,\ldots,13 \}$, we
kept only tweets that have an event specific term and used those tweets
to construct the {{corresponding}} retweet graph, {{denoted by}}
$G_E$. Our analysis focuses on the giant component of the retweet
graph, {{denoted by}} $G_E^0$. In Table~\ref{table:data} we present
important properties of each retweet graph's giant component
{{including}} the number of vertices, number of edges, maximal degree,
and the Twitter name of the superstar corresponding to the maximal
degree. A more detailed description of each event, including the event
specific term, can be found in the \hyperref[app]{Appendix}.

The sizes of the giant components range from 239 to 7365 vertices.
The giant components {{of the retweet graphs corresponding to these
events}} are not trees, but they are very tree-like in that they have
only a few small cycles. In Table~\ref{table:data}, one sees that
for each giant component, the deletion of a small number of edges will
result in an honest tree.

%
\begin{table}
\caption{For each event $E$, we list the number of vertices
[$|V(G_E^0)|$], number of edges [$|E(G_E^0)|$] and maximal degree
[$d_{\mathrm{max}}(G_E^0)$] in the giant component $G_E^0$, along with the
Twitter name of the superstar corresponding to the maximal degree}
\label{table:data}
\begin{tabular*}{\textwidth}{@{\extracolsep{\fill}}lcccc@{}}
\hline
\multicolumn{1}{@{}l}{$\bolds{E}$} & \multicolumn{1}{c}{$\bolds{|V(G_E^0)|}$} &
\multicolumn{1}{c}{$\bolds{|E(G_E^0)|}$} &
\multicolumn{1}{c}{$\bolds{d_{\mathrm{max}}(G_E^0)}$} & \multicolumn{1}{c@{}}{\textbf{Superstar}}\\
\hline
\phantom{0}1 & 7365 & 7620 & 512 & warrenellis\\
\phantom{0}2 & 3995 & 4176 & 362 & anison\\
\phantom{0}3 & 2847 & 2918 & 566 &FIFAWorldCupTM\\
\phantom{0}4&2354& 2414& 657 &taytorswift13\\
\phantom{0}5 &1897& 1929& 256& FIFAcom\\
\phantom{0}6 & 1724 & 1814 & 992 & ladygaga\\
\phantom{0}7& 1659& 2059 & \phantom{0}56& MMFlint\\
\phantom{0}8&1408& 1459 & 269&FIFAWorldCupTM\\
\phantom{0}9&1025& 1045 & 247& FIFAWorldCupTM\\
10& 1024& 1050 & 229& SkyNewsBreak\\
11 &\phantom{0}705& \phantom{0}710 & 113& realmadrid\\
12&\phantom{0}505& \phantom{0}521 & 186& Wimbledon\\
13& \phantom{0}239& \phantom{0}247 & \phantom{0}38& cnnbrk\\
\hline
\end{tabular*}
\end{table}

%
%
\subsection{Maximal degree}
\begin{figure}

\includegraphics{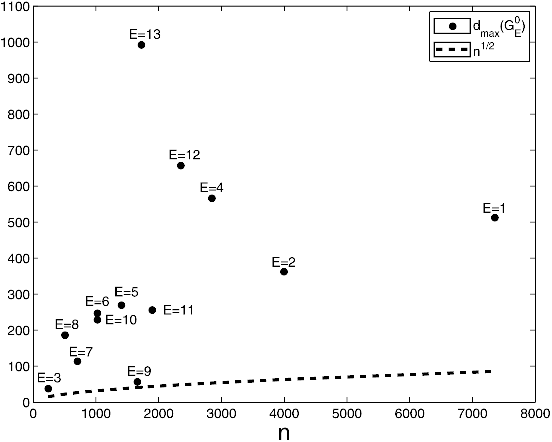}

\caption{Plot of $d_{\mathrm{max}}(G_E^0)$ versus $n=|V(G_E^0)|$ for the
retweet graphs of each event. The events are labeled with the same numbers
as in Table \protect\ref{table:data}. Also shown is a plot of $n^{1/2}$.}
\label{fig:superstar}
\end{figure}

The maximal degree in the retweet graphs is larger than would be
expected under preferential attachment. {{Write $n=|V(G_E^0)|$ for the
number of vertices in the giant component.}} For a \pa graph with $n$
vertices, it is known that the maximal degree scales as $n^{1/2}$.
Figure~\ref{fig:superstar} shows a plot of the maximal degree in
the giant component $d_{\mathrm{max}}(G_E^0)$ and a plot of $n^{1/2}$
versus $n$ for the retweet graphs.
It can be seen from the figure that the sublinear growth
predicted by \pa{{does not capture}}
the superstar effect in these retweet graphs.

%


%
%
\subsection{Estimating $p$ and the degree distribution}

\begin{figure}

\includegraphics{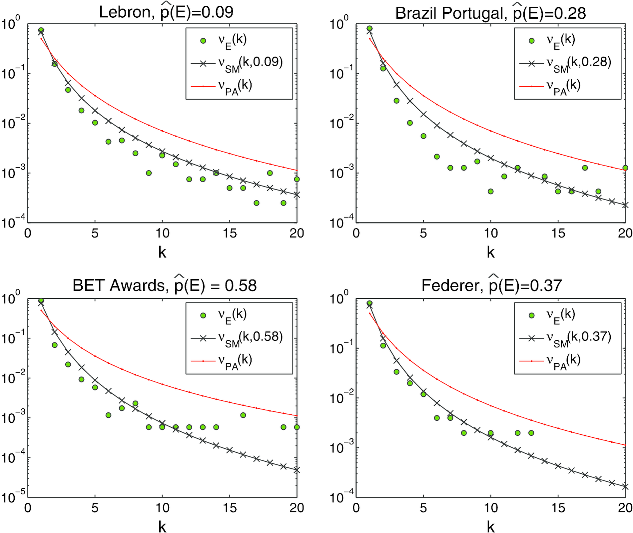}

\caption{Plots of the empirical degree distribution for the giant component
of the retweet graphs [$\nu_E(k)$], and the estimates of the Superstar model
[$\nu_{\mathrm{SM}}(k,\widehat{p}(E))$] and \pa[$\nu_{\mathrm{PA}}(k)$]
for four different events. Each plot is labeled with
the event specific term and $\widehat{p}(E)$.}
\label{fig:degreePMF}
\end{figure}

\begin{figure}

\includegraphics{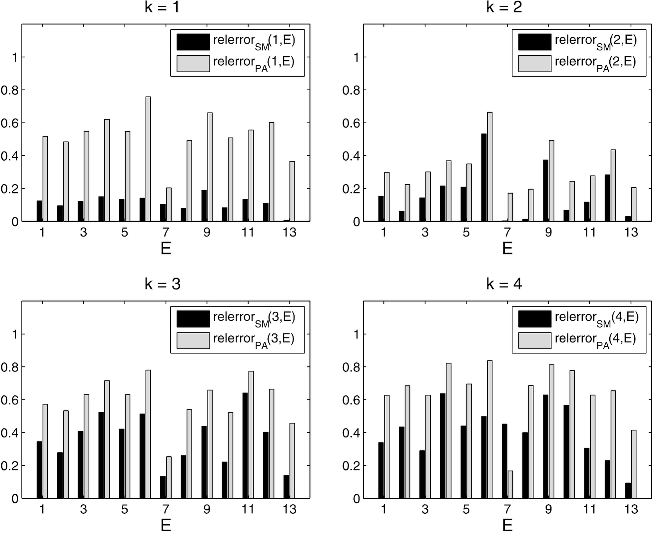}

\caption{Plots of the relative errors of
the degree distribution predictions of \pa and
the Superstar model for 13 retweet graphs. The errors are
plotted for degree $k=1,2,3,4$.}
\label{fig:errorK}
\end{figure}

The {{asymptotic degree distribution of the Superstar model is known
(via Theorem~\ref{thm:DegreeDist})}} once
the superstar parameter $p$ is specified.
We {{were}} interested in seeing,
for each event $E$, {{how well this model predicted the observed degree
distribution}} in $G_E^0$.
For an event $E$ and degree $k\in\{1,2,\ldots\}$,
we define the empirical
degree distribution of the giant component as
\[
\widehat{\nu}_E(k)=\frac{1}{|V(G_E^0)|} \operatorname{card} \bigl
\{
v_j\in V\bigl(G_E^0\bigr)\dvtx\deg
\bigl(v_j, G_E^0\bigr)=k \bigr\}.
\]
To predict the degree distribution using the Superstar model,
we need a value for $p$. We estimate $p$ for each event $E$ as
$\widehat{p}(E)=d_{\mathrm{max}}(G_E^0)/|V(G_E^0)|$.
Using $p=\widehat{p}(E)$ we obtain the Superstar
model degree distribution prediction for each
event $E$ and degree $k$, $\nu_{\mathrm{SM}}(k,\widehat{p})$ from Theorem~\ref
{thm:DegreeDist}.
For comparison, we also compare $\widehat{\nu}_E(k)$ to
the \pa degree distribution $\nu_{\mathrm{PA}}(k)=4 (k(k+1)(k+2) )^{-1}$~\cite{Bollobas:2001:DSS:379831.379835}.
Figure~\ref{fig:degreePMF}
shows the empirical degree distribution for the retweet
graphs of four of the
events, along with the predictions
for the two models.
As can be seen, the Superstar model predictions
seem to qualitatively match the empirical degree
distribution better than preferential attachment.
To obtain a more quantitative comparison of the degree distribution,
we calculate the relative error of these models
for each value of degree $k$. The relative error
for event $E$ and degree $k$ is defined as
$\operatorname{relerror}_{\mathrm{SM}}(k,E) = \llvert\nu_{\mathrm{SM}}(k,\widehat
{p})-\widehat{\nu
}_E(k)\rrvert(\nu_{\mathrm{SM}}(k,\widehat{p}))^{-1}$
for the Superstar model and $\operatorname{relerror}_{\mathrm{PA}}(k,E) =
\llvert\nu
_{\mathrm{PA}}(k)-\widehat{\nu}_E(k)\rrvert(\nu_{\mathrm{PA}}(k))^{-1}$
for preferential attachment.
In Figure~\ref{fig:errorK}, we show the relative errors
for different values of $k$.
As can be seen, the relative error of the Superstar model
is lower than \pa for degrees $k=1,2,3,4$ and for all of the events
with the exception of $k=4$ and $E=7$.
There is a clear connection between the
superstar degree and the degree distribution
in the giant component of these retweet graphs
that is captured well by the Superstar model.

%

%

\section{Analysis of a special two-type branching process}\label{sec:proofs}

{{Let us now start}} the proofs of the main theorems of
Section~\ref
{sec:results}. {{The core of the proof is}} a special two-type
continuous time branching processes together with a surgery operation
that establishes the equivalence between this {{continuous time}}
construction and the {{original}} Superstar model. We start by
describing this construction and {{then}} prove the equivalence
between the two models. 

\subsection{A two-type continuous branching process}
\label{constr:multi-bp}

We now consider a two-type continuous time branching process $\BP(t)$ whose
types we call red and blue. For each fixed $t\geq0$, we shall view
$\BP
(t)$ as {{a}} random tree representing the genealogical {{structure}}
of the population till time $t$. {{This}} includes
parent child relationships of vertices as well as the color of each
vertex. We use $|\BP(t)|$ for the total number of individuals in the
population by time $t$. Every individual survives forever. We shall
also let $ \{\BP(t) \}_{t\geq0}$ be the associated
filtration of the
process. {{Let us now describe the construction.}} At time $t=0$, we
begin with a
single red vertex that we call $v_1$. For any fixed time $0<t< \infty$,
let $V_t$ denote the vertex set of $\BP(t)$. Each vertex
$v\in V_t$ in the branching process gives birth according to a Poisson process
with rate
\[
\lambda(v,t) = c_B(v,t)+1,
\]
where $c_B(v,t)$ is equal to the number of
blue children of vertex $v$ at time $t$. Also let $c_R(v,t)$ denote the
number of red children of vertex {{$v$}} by time $t$. At the moment
of a new birth, {{this new}} vertex
is colored red with probability $p$ and colored blue with
probability $q=1-p$. Finally, for $n\geq1$, {{define}} the
stopping times
%
%
\begin{equation}
\label{eqn:stop-time-def} \tau_n = \inf\bigl\{t\dvtx\bigl|\BP(t)\bigr| = n
\bigr\}.
\end{equation}
Since the counting process $|\BP(t)|$ is a nonhomogenous Poisson
process with
a rate that is always greater than or equal to one,
the stopping times $\tau_n$ are almost surely finite. {{This completes
the construction of the branching process.}}

%

\subsection{Equivalence between the branching process and the
Superstar model}
\label{sec:equiv}

%
%
Before diving into properties of our two-type branching process
constructed as above, let us show how the Superstar model can be
obtained from the above branching process via a \emph{surgery} operation.
We start with an informal description of the connection between the
Superstar model and the branching
process $\BP(\cdot)$. {{To describe this connection, we introduce a
new vertex $v_0$ namely the superstar vertex to the system. Recall that
$v_1$ was the root (the initial progenitor) of the branching process
$\BP(\cdot)$.}} We connect vertex $v_1$, to the superstar
$v_0$ [{{$v_0$ played no role in the evolution of $\BP(\cdot)$}}].
This forms the Superstar model $G_2$ on $2$ vertices. All the red
vertices in the process $\BP(\cdot)$ correspond to the neighbors of the
superstar $v_0$. The true degree of a (nonsuperstar) vertex in
$G_{n+1}$ is one plus the number of its blue children
in $\BP(\tau_n)$, where the additional factor of one comes from the
edge connecting this vertex to its parent.
By elementary properties of the exponential distribution, the dynamics
of $\BP(\cdot)$ imply that each new vertex that is born is red
(connected to the superstar $v_0$) with probability $p$, else with
probability $q=1-p$ is blue and
connected to one of the remaining {{extant (nonsuperstar)}}
vertices with probability proportional to the current degree of that
vertex, thus
increasing the degree of this chosen vertex by one.
These dynamics are the same as the Superstar model.

Formally, our surgery will take the tree $\BP(\tau_{n})$ and modify it
to get an $(n+1)$-vertex tree $\mathS_n$ that has the same distribution
as the Superstar model $G_{n+1}$.
From this, we will be able to read off the probabilistic properties of
the superstar tree $G_{n+1}$.

\begin{figure}

\includegraphics{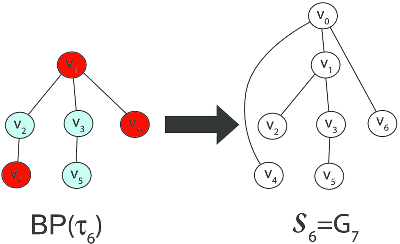}

\caption{The surgery passing from $\BP(\tau_{n})$ to $\mathS_{n+1}$ and
$G_{n+1}$ for $n=6$.}
\label{fig:surgery}
\end{figure}

We label the vertices of $\BP(\tau_{n})$ as $\{v_1, v_2,\ldots,
v_n\}$
in order of their birth. {{Now add a new vertex $v_0$ to this set to
give us the vertex set of the tree $\cS_{n}$}}. One can
anticipate that $v_0$ will be our superstar. Next, we define the edge
set for $\mathS_n$. To do this, we
take each red vertex $v$ in $\BP(\tau_{n})$, remove the
edge connecting $v$ to its parent (if it has one) and then we
create a new edge between $v$ and $v_0$. To complete the construction
of $\mathS_n$,
it only remains to ignore the color of the vertices. An illustration of
this surgery for $n=6$ is given
in Figure~\ref{fig:surgery}.

%

%
\begin{Proposition}[(Equivalence from surgery operation)]
\label{prop:equiv}
The sequence of trees $ \{\mathS_n\dvtx n\geq1 \}$ has the
same distribution as the Superstar model $ \{G_{n+1}\dvtx\break  n\geq
1 \}$.
\end{Proposition}

\begin{pf} Think of $\mathS_n$
as being rooted at $v_0$ so that every vertex except $v_0$ in $\mathS_n$
has a unique parent. The parent of all the red individuals is the
superstar $v_0$ while the parents of all of the other blue individuals
are unchanged from $\BP(\tau_n)$.

The induction hypothesis will be that $\mathS_n$ has the same
distribution as $G_{n+1}$ and the degree of each nonsuperstar vertex in
$\mathS_n$ is the number of blue children it possesses plus one for the
edge connecting the vertex to its parent in $\mathS_n$. Condition on
$\BP(\tau_n)$ and fix $v\in\BP(\tau_n)$. By the property of the
exponential distribution, the probability that the next vertex born
into the system is born to vertex $v$ is given
\[
\frac{\lambda(v,\tau_n)}{\sum_{u\in\BP(\tau_n)}\lambda({{u}},
\tau
_n)} = \frac{c_B(v, \tau_n)+1}{\sum_{u\in\BP(\tau_n)}c_B({{u}},
\tau_n)+1}.
\]
Thus a new vertex {{$v_{n+1}$}} attaches to vertex $v$ with
probability proportional
to the present degree of $v$ in $\mathS_n$. Further, with probability
$p$, this vertex
is colored red, whence by the surgery operation, the edge to
{{$v_{n+1}$}} is deleted and this new vertex is connected to the
superstar $v_{0}$. In this case the degree of $v$ in $\mathS_n$ is
unchanged. With probability $1-p$ this new vertex is colored blue,
whence the surgery operation does not disturb this vertex so that the
degree of vertex $v$ is increased by one. These are exactly the
dynamics of $G_{n+2}$ conditional on $G_{n+1}$. By induction the result
follows.
\end{pf}

For the rest of the paper, we shall assume $G_{n+1}$ is constructed
through this surgery process {{from $\BP(\tau_n)$}} and suppress
$\mathS_n$.

\subsection{Elementary properties of the branching process}
The previous section {{set up}} an equivalence between the Superstar
model and the two type continuous time branching process. The aim of
this section is to prove properties of this two type branching process.
Section~\ref{sec:comp-proof} uses these results to complete the proof
of the main results for the Superstar model.

For $t\geq0$, write $R(t)$ and $B(t)$ for the total number of red and
blue vertices, respectively, in $\BP(t)$. By construction of the
process $ \{\BP(t)\dvtx t\geq0 \}$, every new vertex is
independently
colored red with probability $p$ and blue with probability $1-p$. In
particular, the number of blue vertices $B(t)$ is just a time change of
a random walk with Bernoulli($1-p$) increments. Thus, by the strong law
of large numbers
%
%
\begin{equation}
\label{eqn:ss-b-convg} \frac{B(t)}{|\BP(t)|} \stackrel
{\mathrm{a.s.}}{\longrightarrow}1-p\qquad \mbox{as } t\to\infty.
\end{equation}
Before moving onto an analysis of the branching process, we introduce
the Yule process.
%
%
\begin{Definition}[(Rate $a$ Yule process)]\label{def:yule}
{{Fix $a> 0$. A}} rate $a$ Yule process is defined as a pure birth
process $\Yu_a(\cdot)$ that starts with a single individual $\Yu
_a(0)=1$ {{and}} with the rate of creating new individuals
proportional to the number of present individuals in the population, namely
\[
\prob\bigl(\Yu_a(t+dt) - \Yu_a(t) = 1 | \bigl\{
\Yu_a(s)\dvtx0\leq s\leq t \bigr\}\bigr) = a \Yu_a(t)
\,dt.
\]
\end{Definition}
The Yule process is a well-studied {{probabilistic}} object. {{The}}
next lemma collects some of its standard properties. In
particular, part (a) follows from \cite{norris1998markov}, Section~2.5,
while (b) follows from \cite{ref:an}, Theorem~1, III.7.
%
%
\begin{Lemma}[(Yule process)]\label{lem:yule-proc}
\textup{(a)} For each $t> 0$, the random variable $\Yu_a(t)$ has a geometric
distribution with parameter $e^{-at}$, that is,
\[
\prob\bigl(\Yu_a(t) = k \bigr) = e^{-at}
\bigl(1-e^{-at}\bigr)^{k-1},\qquad k\geq1.
\]
\textup{(b)} The process $ (e^{-at} \Yu_a(t)\dvtx0\leq t< \infty
)$ is an
$\mathbb{L}^2$ bounded martingale with respect to the natural
filtration and $e^{-at} \Yu_a(t) \stackrel{a.s.}{\longrightarrow
}W^\prime$, where $W^\prime$
has an
exponential distribution with mean one.
\end{Lemma}
Now define the process
\[
M(t)= e^{-(2-p)t} \bigl(\bigl|\BP(t)\bigr|+B(t) \bigr), \qquad t\geq0.
\]
Note that $M(0) =1$.
%
%
\begin{Proposition}[{[Asymptotics for $\BP(t)$]}]\label{pro:conts-mp}
The process $ (M(t)\dvtx t\geq0 )$ is a positive $\mathbb{L}^2$
bounded martingale with respect to the natural filtration $ \{\BP
(t)\dvtx{t\geq0} \}$, and thus converges almost surely and in
$\mathbb
{L}^2$ to
a random variable $W^*$ with $\expec(W^*) =1$. The random variable
$W^*$ is positive with probability one. {{Further, one has}}
%
%
\begin{equation}
\label{eq:FFLimitAS} \lim_{t \rightarrow\infty} {e^{-(2-p)t}\bigl|\BP
(t)}\bigr| =
\frac{W^*}{2-p} \qquad\mbox{with probability one}.
\end{equation}
\end{Proposition}
%

\begin{pf} We write $Z(t) = |\BP(t)|$ and $Y(t) = Z(t)+B(t)$ so that
$M(t) = e^{-(2-p)t} Y(t)$ and we let $dM(t)=M(t+dt)-M(t)$. We then have
%
%
\begin{equation}
\label{eqn:mt-mart-decomp} dM(t) = e^{-(2-p)t} \,dY(t) -
(2-p)e^{-(2-p)t} Y(t) \,dt.
\end{equation}
The processes $Z(t), B(t)$ are all counting processes. For such
processes, we shall use the infinitesimal shorthand $\expec(dZ(t)|\BP
(t)) = a(t) \,dt$ to denote the fact that $Z(t) - \int_0^t a(s)\,ds$
is a
local martingale.

Now the counting process $Z(t)=|\BP(t)|$ evolves by jumps of size one with
%
%
\begin{equation}
\label{eqn:dzt-evol-eqn} \prob\bigl( dZ(t)=1|\BP(t) \bigr) =
\biggl(\sum
_{v\in\BP
(t)} \bigl(c_B(v,t)+1\bigr) \biggr) \,dt,
\end{equation}
where $c_B(v,t)$ always denotes the number of blue children of vertex
$v$ at time $t$. The number of blue vertices can be written as $B(t) =
\sum_{v\in\BP(t)} c_B(v,t)$ since every blue vertex is an offspring of
a unique vertex in $\BP(t)$. {{Using \eqref{eqn:dzt-evol-eqn} results in}}
\[
\expec\bigl(dZ(t)|\BP(t)\bigr) = \bigl(Z(t)+ B(t)\bigr) \,dt.
\]
Since $B(t)\leq Z(t)$, we see that the rate of producing new
individuals is bounded by $2|\BP(t)|$. Thus, the process $|\BP(t)|$ can
be stochastically bounded by a Yule process with $a=2$. This implies by
Lemma~\ref{lem:yule-proc} that for all $t\geq0$ we have $\expec(|\BP
(t)|^2) < \infty$.

Let us now analyze the process $B(t)$. This process increases by one
when the new vertex born into $\BP(\cdot)$ is colored blue that happens
with probability $1-p$. Thus, we get
\[
\expec\bigl(dB(t)|\BP(t)\bigr) = (1-p) \bigl(Z(t)+B(t)\bigr) \,dt.
\]
Combining the last two equation gives us
\[
\expec\bigl(dY(t)|\BP(t)\bigr) = (2-p) Y(t) \,dt.
\]
{{Using \eqref{eqn:mt-mart-decomp} now gives that}} $\expec
(dM(t)|\BP
(t)) = 0$. This completes the proof that $M(\cdot)$ is a martingale.

Next, we check that $M(\cdot)$ is an $\mathbb{L}^2$ bounded
martingale. Since $Y^2(t+dt)$ can take values $(Y(t)+1)^2$ or
$(Y(t)+2)^2$ at rate $pY(t)$ and $(1-p)Y(t)$, respectively, we have
\[
\expec\bigl(d\bigl(M^2(t)\bigr)|\BP(t)\bigr) = (4-3p)
e^{-(2-p)t} M(t)\,dt.
\]
{{Thus,}} the process $U(t)$ defined by
\[
U(t) = M^2(t) - (4-3p)\int_0^t
e^{-(2-p)s} M(s)\,ds,
\]
is a martingale. Taking expectations and noting that since
$M(\cdot)$ is a martingale, with $M(0)=1$ thus $\expec(M(s)) = 1$ for
all $s$, we get
\[
\expec\bigl(M^2(t)\bigr) = 1+ (4-3p)\int_0^t
e^{-(2-p)s}\,ds \leq1 + \frac
{4-3p}{2-p}.
\]
This $\mathbb{L}^2$ boundedness implies
that there exists a random variable $W^*$ such that
\[
e^{-(2-p)t} \bigl(\bigl|\BP(t)\bigr|+B(t)\bigr) \stackrel{\mathrm{a.s.}, \mathbb
{L}^2} {\longrightarrow} W^*.
\]
{{Using}} \eqref{eqn:ss-b-convg} shows that $e^{-(2-p)t} |\BP(t)|
\to W^*/(2-p)$. To ease notation, write
\[
W:= \frac{W^*}{(2-p)}.
\]
To complete the proof of the proposition we need to show that
$W$ is strictly positive. First, note that by
$\mathbb{L}^2$ convergence, $\expec(W^*) =1$. So in particular
$\prob(W = 0)=r <1$. Let $\zeta_1< \zeta_2< \cdots$
be the times of birth of children (blue or red) of the root
vertex $v_1$ and write $\BP_i(\cdot)$ for the
subtree consisting of the $i$th child and its descendants. Then
\[
e^{-(2-p)t}\bigl|\BP(t)\bigr| = \sum_{i=1}^\infty
e^{-(2-p)\zeta_i} \bigl[e^{-(2-p)(t-\zeta_i)}\bigl|\BP_i(t-
\zeta_i)\bigr| \bigr] \ind\{\zeta_i\leq t \} +
e^{-(2-p)t}.
\]
Thus, as $t\rightarrow\infty$, for any fixed $K\geq1$, we have
\[
W \geq_{\mathrm{st}}\sum_{i=1}^K
e^{-(2-p)\zeta_i} W_i,
\]
where $ \{W_i \}_{i\geq1}$ are independent and identically
distributed with the same distribution as $W$
(independent of $ \{\zeta_i \}_{i\geq1}$) and $\geq_{\mathrm{st}}$ denotes
stochastic domination. This independence gives us
\[
\prob(W = 0) \leq\prob(W_i = 0\ \forall1\leq i\leq K) =
r^K.
\]
{{Letting}} $K\to\infty$ that $\prob(W=0)=0$.
%
\end{pf}

Before ending this section, we derive some elementary properties of the
offspring of an individual in $\BP(\cdot)$. Let $\sigma_v$ be the time
of birth of vertex $v$ in $\BP(\cdot)$. Recall that $c_B(v, \sigma
_v+s)$ and $c_R(v, \sigma_v+s)$ denote the number of blue and red
children, respectively, of this vertex $s$ units of time after the birth
of $v$. {{Since the distribution of the point process representing
offspring of each vertex is the same}}, these random variables have the
same distribution irrespective of the choice of the vertex $v$. Define
the process
\[
M^*(t):= c_R(v,{{\sigma_v+t}}) - {{p}}\int
_0^t \bigl(c_B(v,
\sigma_v+ s)+1\bigr) \,ds,\qquad t\geq0.
\]
%
%
\begin{Lemma}[({{Offspring point process: distributional
properties}})]\label{lemma:offspring}
\begin{longlist}[(a)]
\item[(a)] Conditional on $\BP(\sigma_v)$ we have
\[
\bigl(c_B(v, \sigma_v+t)\dvtx t\geq0 \bigr)
\stackrel{d} {=} \bigl(\Yu_{1-p}(t) - 1\dvtx t\geq0 \bigr),
\]
and thus one has
\[
\expec\bigl(c_B(v,\sigma_v+ t)\bigr) =
e^{(1-p)t} - 1,\qquad t\geq0.
\]

\item[(b)] The process $ (M^*(t)\dvtx t\geq0 )$ is a martingale with
respect to
the filtration $ \{\BP(\sigma_v+t)\dvtx t\geq0 \}$ and
one has
\[
\expec\bigl(c_R(v, \sigma_v+t)\bigr) =
\frac{p}{1-p}\bigl(e^{(1-p)t} -1\bigr), \qquad t\geq0.
\]
\end{longlist}
\end{Lemma}
\begin{pf} Part (a) is obvious from construction. To prove (b),
note that
\[
\expec\bigl(d c_R(v, \sigma_v+ t)|\BP(t+
\sigma_v)\bigr) = p\bigl(c_B(v,\sigma_v+
t)+1\bigr) \,dt,
\]
since vertex $v$ creates a new child at rate $c_B(v,\sigma_v+t)+1$
which is then marked red with probability $p$.
\end{pf}

\subsection{Convergence for blue children proportions}
The equivalence between $\BP(\cdot)$ and the Superstar model
{{described in Section~\ref{sec:equiv}}} will
{{imply}} that the number of vertices with degree $k+1$ in $G_{n+1}$
is the same as the number of vertices in $\BP(\tau_n)$ with exactly
$k$ blue children. We will need general results on the asymptotics
of such counts for the process $\BP(t)$ as $t\to\infty$. {{Using the
equivalence created by the surgery operation, one can then transfer
these results to asymptotics for the degree distribution of the
original Superstar model.}} {{Now}} recall the random variable $W^*$
obtained as the martingale limit obtained in Proposition~\ref
{pro:conts-mp}. Define
$p_{\geq k}(\infty)$ as
%
%
\begin{equation}
\label{eqn:p-geqk-defn} p_{\geq k}(\infty)=k!\prod_{i=1}^k
\biggl(i+\frac{2-p}{1-p} \biggr)^{-1}.
\end{equation}
%
%
\begin{Theorem}\label{thm:degree-count-conv}
Fix $k\geq1$ and let $Z_{\geq k}(t)$ denote the number of vertices in
$\BP(t)$ that have at least $k$ blue children. Then
\[
e^{-(2-p)t} Z_{\geq k}(t) \stackrel{a.s.}{\longrightarrow}p_{\geq
k}(\infty)
\frac{W^*}{2-p}
\]
as $t\to\infty$.
%
\end{Theorem}

\begin{pf} The proof uses a variant of the reproduction martingale
technique developed in \cite{nerman1981convergence} and it is framed in
two steps:

\begin{longlist}[(a)]
\item[(a)] Proving convergence of expectations of the desired quantities to
the expectations of the asserted limits. This is proved in Section~\ref
{sec:convg-expec}.
\item[(b)] Bootstrapping this convergence to almost sure convergence using
laws of large numbers. This is proved in Section~\ref{sec:as-convg}.
\end{longlist}

We {{start with}} some notation {{required}} to carry out this
program. For a vertex $v$, write
\[
\bolds{\zeta}^v = \bigl(\bigl(\xi_i^v,
\mathcal{C}_i^v\bigr)\dvtx i\geq1\bigr),
\]
for the {{point process representing offspring (times of birth and
types) of this vertex~$v$. More precisely}} here $\xi_i^v$ denotes the
time of birth of the $i$th {{offspring}} of vertex $v$ \emph{after}
{{the birth of vertex $v$}} into the branching process $ \{\BP
(t)\dvtx t\geq0 \}$ while $\mathcal{C}_i^v$ denotes the color of this child
(red or
blue). {{Thus, the $i$th offspring of vertex $v$ is born into $\BP$ at
time $\sigma_v + \xi_i^v$.}}
Write $\bolds{\xi}^v =({{\xi_i^v}}\dvtx i\geq1 )$ for the process
that just
keeps track of times of birth of these offspring for vertex $v$. {{Note
that the point processes $\bolds{\zeta}^v$ and $\bolds{\xi}^v$ have
the same
distribution across vertices $v$.}} We shall use $\bolds{\zeta}:=
\bolds{\zeta}
^{v_1}$ and $\bolds{\xi}:= \bolds{\xi}^{v_1}$ to denote {{a}}
generic point
process with the above distributions. We shall view $\bolds{\xi}$ as a
counting measure on $(\mathbb{R}_+, \mathcal{B}(\mathbb{R}_+))$. For
$A\in\mathcal{B}(\mathbb{R}_+)$, write
$\bolds{\xi}(A)$ for the number of points in the set $A$. {{Define the
corresponding intensity measure}} $\mu$ {{by}}
\[
\mu(A):= \mathbb{E}\bigl(\bolds{\xi}(A)\bigr), \qquad A\in\mathcal
{B}(\mathbb{R}_+).
\]
We start with a simple lemma that has notable consequences.

%
\begin{Lemma}[(Renewal measure)]\label{lem:renwal-measure}
For $\alpha=2-p$, we have
\[
\int_0^\infty e^{-\alpha t} \mu(dt) = 1.
\]
The measure defined by setting $\mu_\alpha:= e^{-\alpha t} \mu(dt)$
is a probability measure and this measure has expectation
$\int_0^\infty t \mu_\alpha(dt) = 1 $.
\end{Lemma}

\begin{pf} As in Lemma~\ref{lemma:offspring},
let $c_B(v_1,t)$ and $ c_R(v_1, t)$ denote the number
of red and blue children, respectively, of vertex $v_1$ by time $t$
{{(note that $\sigma_{v_1} = 0$ )}}.
Then {{by definition, the intensity measure $\mu$ satisfies}} $\mu
([0,t]) = \expec(c_R(v_1, t) +c_B(v_1, t))$.
Further by Fubini's theorem,
\[
\int_0^\infty e^{-\alpha t} \mu(dt) = \alpha
\int_0^\infty e^{-\alpha
t}\mu[0,t] \,dt.
\]
Using the expressions for $\expec(c_B(v_1,t))$ and $ \expec
(c_R(v_1,t))$ from Lemma~\ref{lemma:offspring} completes the proof. The
second assertion regarding the expectation follows similarly.
\end{pf}

\subsubsection{Convergence of expectations}\label{sec:convg-expec}
The first step in the proof of Theorem~\ref{thm:degree-count-conv} is
convergence of expectations. This follows using standard renewal
theory. We setup notation that allows us to use the linearity of
expectations to derive a renewal equation. We start with the definition
of a \textit{characteristic} \cite{jagers-book,jagers-nerman-comp}
that we
use to count the number of vertices in the branching process with some
fixed property. For each vertex $v\in\BP(\infty)$, let $ \{
\phi^v(s)\dvtx s\geq0 \}$ be an independent and identically distributed
nonnegative stochastic process, with $\phi^v(s)$ measurable with
respect to $ \{(\xi_i^v,\mathcal{C}_i^v)\dvtx\xi_i^v\leq s \}$.
Thus, the
value of
the stochastic process at time $s$ {{namely}} $\phi^v(s)$ is
determined by the set offspring of vertex $v$ born before the age $s$
of this vertex~$v$.

The value $\phi^v(s)$ is referred to as the score of vertex $v$ at age
$s$ \cite{jagers-book}, Section~6.9. We write $\phi:=\phi^{v_1}$ to
denote the process corresponding to the root when we would like to
refer to a generic such process. Throughout we shall assume that $\phi
(\cdot)$ is bounded and nonnegative, namely for some constant $C<
\infty$,
\[
\phi(s)\geq0, \qquad\phi(s) < C\qquad \mbox{for all } s\geq0.
\]
Define
\[
Z_\phi(t) = \sum_{v\in\BP(t)}
\phi^v(t-\sigma_v),\qquad  t\geq0
\]
for the branching process $\BP(\cdot)$ counted according to
characteristic $\phi$. The main examples of interest are:
\begin{longlist}[(a)]
\item[(a)] \textit{Total size}: $\phi(s) = 1$ for all $s\geq0$. This results
in $Z_\phi(t) =|\BP(t)|$, the total size of the branching process by
time $t$.

\item[(b)] \textit{Degree}: $\phi(s) = \ind\{k \mbox{ or more blue
children at age $s$} \}$ gives $Z_\phi(t) =\break  Z_{\geq k}(t)$, the
number of vertices
in $\BP(t)$ with $k$ or more blue children.
\end{longlist}

Now fix an \textit{arbitrary} bounded characteristic $\phi$. For fixed
time $t> 0$, conditioning on the offspring process $\bolds{\zeta}
:=\bolds{\zeta}^{v_1}$ of vertex ${v_1}$, the branching process counted
according to this characteristic satisfies the recursion
%
%
\begin{equation}
\label{eqn:charac-recursion} Z_\phi(t) = \phi^{v_1}(t) + \sum
_{\xi_i^{v_1}\leq t } Z_\phi^{\sss
(i)}\bigl(t-
\xi_i^{v_1}\bigr),
\end{equation}
where $Z_\phi^{\sss(i)}(\cdot) \stackrel{d}{=} Z_\phi(\cdot)$ and are
independent {{for $i\geq1$ and correspond to the contribution of the
descendants of the $i$th child of vertex $v_1$.}} Taking
expectations and {{defining the function $m_{\phi}(\cdot)$ by}}
$m_\phi(t) := \expec(Z_\phi(t))$, {{this function satisfies}} the
renewal equation
\[
m_\phi(t) = \expec\bigl(\phi(t)\bigr) + \int_0^t
m_\phi(t-s) \mu(ds).
\]
Define
\[
\tilde{m}_\phi(t):= e^{-\alpha t} m_\phi(t),\qquad t\geq0.
\]
Lemma~\ref{lem:renwal-measure} and {{standard}} renewal theory
(\cite{jagers-book}, Theorem~5.2.8)
now imply the next result.
%
%
\begin{Proposition}\label{pro:renewal-expec}
For arbitrary bounded characteristics, writing $\alpha= (2-p)$ we have
\[
\lim_{t\to\infty} \tilde{m}_\phi(t) {{=}} \int
_0^\infty e^{-\alpha s} \expec\bigl(\phi(s)
\bigr) \,ds:= \tilde{m}_\phi(\infty).
\]
\end{Proposition}
Applying this to the two examples which count the size of the branching
process and number of vertices with at least $k$ blue children, we get
the following result.
%
%
\begin{Corollary}
\label{cor:pk-expression}
Taking the two characteristics of interest one gets for $\phi(t) = 1$
\[
e^{-\alpha t}\expec\bigl(\bigl|\BP(t)\bigr|\bigr) \to\frac{1}{\alpha}\qquad
\mbox{as } t
\to\infty
\]
and for $\phi(t) = \ind\{k \mbox{ or more blue children at
time $t$} \}$
\[
e^{-\alpha t}\expec\bigl(Z_{\geq k}(t)\bigr) \to\frac{p_{\geq
k}(\infty)}{\alpha}\qquad
\mbox{as } t\to\infty,
\]
with $p_{\geq k}(\infty)$ as in \eqref{eqn:p-geqk-defn}.
\end{Corollary}
\begin{pf} The first assertion in the corollary is obvious
{{[corresponding to the case $\phi(\cdot) \equiv1 $]}}.
To prove the second assertion regarding the number of blue
vertices, observe that the limit constant in
Proposition~\ref{pro:renewal-expec} can be written as
\begin{eqnarray*}
&&\frac{1}{\alpha}\int_0^\infty\alpha
e^{-\alpha s} \expec\bigl(\ind\{\mbox{root } {{v_1}} \mbox{ has
} k \mbox{ or more blue children at age $s$} \}\bigr)\,ds\\
&&\qquad = \frac
{1}{\alpha}\prob
\bigl(c_B(v_1,T)\geq k\bigr),
\end{eqnarray*}
where $T$ is an exponential random variable with mean $\alpha^{-1}$
that is independent of the counting process of the number blue
offspring $c_B(v_1,\cdot)$. {{Further, by Lemma~\ref
{lemma:offspring}(a),}}
\[
c_B(v_1,\cdot) \stackrel{d} {=} \Yu_{1-p}(
\cdot)-1,
\]
where $\Yu_{1-p}(\cdot)$ is rate $1-p$ Yule process. The interarrival
times $X_i$ between blue children
$i$ and $i+1$ are independent exponential
random variables with mean $(1-p)^{-1}(i+1)^{-1}$, independent of $T$.\vspace*{1pt}
In particular
$\prob(c_B(v_1, T)\geq k) = \prob(T> \sum_{j=0}^{k-1} X_j)$.
Conditioning on the value of $\sum_{j=0}^{k-1} X_j$ and using tail
probabilities for the exponential distribution shows that
\[
\prob\Biggl(T> \sum_{j=0}^{k-1}
X_j\Biggr)= \expec\Biggl(\exp\Biggl(-\alpha\sum
_{j=0}^{k-1} X_j \Biggr) \Biggr) = \prod
_{j=0}^{k-1} \expec\bigl(\exp(-\alpha
X_j) \bigr).
\]
Using the Laplace transform of the exponential distribution, one can
check that the last expression equals $p_{\geq k}(\infty)$.
\end{pf}

\subsubsection{Almost sure convergence}\label{sec:as-convg}
The aim of this section is to strengthen the convergence of
expectations to almost sure convergence. A key role is played by a
reproduction martingale, a close relative of the martingale used in
\cite{nerman1981convergence} to analyze single type branching processes
as well as in \cite{kingman-age-depen} to analyze times of first birth
in generations. Let $v_1, v_2, v_3, \ldots$ denote the vertices of
$\BP
(\cdot)$ listed in the order of their birth times and let $\sigma
_{v_i}$ denote the time at which vertex $v_i$ is born into the
branching process $\BP(\cdot)$. Note that $\sigma_{v_1} = 0$. Recall
that $\bolds{\xi}^{v_i} = (\xi_{1}^{v_i}, \xi_{2}^{v_i}, \ldots)$
denotes the
offspring point process of $v_i$, namely the first offspring of $v_i$
is born at time $\sigma_{v_i} + \xi_{1}^{v_i} $, the second offspring
of $v_i$ is born at time $\sigma_{v_i} +{{\xi_{2}^{v_i}}} $ and so on.
To ease notation, we shall write $\bolds{\zeta}^{\sss(i)}:=\bolds
{\zeta}^{v_i}$ and
$\bolds{\xi}^{\sss(i)}:= \bolds{\xi}^{v_i}$. Viewing $\bolds{\xi
}^{\sss(i)}$
as a random counting measure on $\mathbb{R}_+$ and writing $\alpha=
2-p$, we have
\[
\xi_\alpha^{\sss(i)} := \sum_{j=1}^\infty
\exp\bigl({-\alpha\xi_{j}^{v_i}} \bigr) = \int
_0^\infty e^{-\alpha t} \bolds{\xi}^{\sss(i)}(dt).
\]
For $m\geq1$, let $\tilde{\mathcal{F}}_m$ be the sigma-algebra
generated by vertices $ \{v_1,\ldots, v_m \} $ and their
offspring processes, namely
\[
\tilde{\mathcal{F}}_m:=\sigma\bigl( \bigl\{\bolds{\zeta}^{\sss
(i)}\dvtx1\leq
i\leq m \bigr\}\bigr).
\]
For $m=0$, let $\tilde{\mathcal{F}}_0$ be the trivial sigma-field.
Now define $\tilde{R}_0 = 1$ and for $m\geq0$ define
\[
\tilde{R}_{m+1} := \tilde{R}_m + e^{-\alpha\sigma_{v_{m+1}}}\bigl
(\xi
_{\alpha}^{\sss(m+1)}-1\bigr).
\]
Let $\Gamma_m$ be the set of the first $m$ individuals born and
\textit{all}
of their offspring. One can check that
%
%
\begin{equation}
\label{eqn:rm-expand} {{\tilde{R}_{m} = \sum_{v\in\Gamma_m}
e^{-\alpha\sigma_v} - \sum_{j=1}^m
e^{-\alpha\sigma_{v_j}}.}}
\end{equation}
Thus, $\tilde{R}_m$ is a weighted sum of children of the first $m$
individuals with weight $e^{-\alpha\sigma_x}$ for vertex $x$, the
individuals $v_1, v_2,\ldots, v_m$ being excluded.
In particular, $\tilde{R}_m > 0$ for all $m$.
The next lemma shows that the sequence {{$(\tilde{R}_m\dvtx{m\geq0})$}}
is much more.
%
%
\begin{Proposition}[(Reproduction martingale)]\label{prop:repr-mart}
The sequence $ (\tilde{R}_m\dvtx{m\geq0} )$ is a nonnegative
$\mathbb{L}^2$
bounded martingale with respect to the filtration $ \{\tilde{\mathcal{F}
}_m\dvtx\break {m\geq0} \}$. Thus, there exists a random variable
$R_\infty$ with
$\expec(R_\infty)=1$ such that $\tilde{R}_m \to R_\infty$ almost surely
and in $\mathbb{L}^2$.
\end{Proposition}
\begin{pf}By the choice of $\alpha= 2-p$ in Lemma~\ref
{lem:renwal-measure} for $i\geq1$, we have $\expec(\xi_\alpha^{\sss
(i)}) = \int_0^\infty e^{-\alpha t} \mu(dt) =1$. Further, $\sigma
_{v_{m+1}}$ is $\tilde{\mathcal{F}}_m$ measurable while $\xi_\alpha
^{\sss
(m+1)}$ is independent of $\tilde{\mathcal{F}}_m$. This implies
\[
\expec(\tilde{R}_{m+1}-\tilde{R}_m|\tilde{
\mathcal{F}}_m) = e^{-\alpha
\sigma
_{v_{m+1}}} \expec\bigl(\xi_\alpha^{\sss(m+1)}
-1\bigr) =0.
\]
By the orthogonality of the increments of the martingale $R_m$, we see that
\[
\expec\bigl((\tilde{R}_m-1)^2\bigr) \leq\expec\bigl(
\bigl[\xi_\alpha^{\sss(i)}\bigr]^2\bigr) \expec\Biggl(
\sum_{i=1}^m e^{-2\alpha\sigma_{v_i}} \Biggr).
\]
Thus, to check $\mathbb{L}^2$ boundedness it is enough to
check that the right-hand side is bounded. The following
lemma bounds the right-hand side of the above equation in two steps and
completes the proof.
\end{pf}
%
%
\begin{Lemma}
\textup{(a)} Let $\xi_\alpha:= \xi^{v_1}_\alpha$ and assume $0< p< 1$. Then
$\expec([\xi_\alpha]^2) < \infty$.

\textup{(b)} For any $m$, $\expec(\sum_{i=1}^m e^{-2\alpha\sigma_{v_i}})
\leq
1+\alpha^{-1}$.
\end{Lemma}
\begin{pf} 
To prove (a), we observe that
$\xi_\alpha= \int_0^\infty\alpha e^{-\alpha t} \bolds{\xi}[0,t] \,
dt $ where
$\bolds{\xi}$ is the point process encoding times of birth of
offspring of
$v_1$. Thus, by Jensen's inequality with the probability measure
$\alpha e^{-\alpha t} \,dt$ we have
\[
[\xi_\alpha]^2\leq\int_0^\infty
\alpha e^{-\alpha t} \bigl[\xi[0,t]\bigr]^2 \,dt.
\]
Let $T$ be an exponential random variable with mean $\alpha^{-1}$
independent of $\xi$. Thus, it is enough to show
$\expec([\xi[0,T]]^2) < \infty$. Note that $\xi[0,T] =
c_R(v_1,T)+c_B(v_1,T)$,
that is, the number of red and blue vertices born to $v_1$ by the
random time $T$.
Thus, it is enough to show $\expec(c_R^2(v_1,T))$
and $\expec(c_B^2(v_1,T)) < \infty$. Conditioning on
$T=t$ {{first note by using}} Lemma~\ref{lem:yule-proc} that
for fixed $t$, $\expec(c_B^2(v_1,t)) \leq C e^{2(1-p)t}$ {{where $C<
\infty$ is a constant independent of $t$.}}
{{Further again using Lemma~\ref{lem:yule-proc}, for any fixed $t$,
conditional on $c_B(v_1, t)$, $c_R(v_1,t)$ is stochastically dominated
by a Poisson random variable with rate $t c_B(v_1,t)$.}} Noting that
$\alpha=2-p$,
we get
\[
\expec\bigl(\bigl[\xi[0,T]\bigr]^2\bigr) \leq C^\prime\int
_0^\infty e^{-(2-p)t} \bigl(e^{2(1-p)t}+
t^2 e^{2(1-p)t} \bigr) \,dt < \infty,
\]
{{for some constant $C^\prime<\infty$. This completes the proof of (a).}}

To prove (b), let $S(t) = \sum_{v\in\BP(t)} e^{-2\alpha\sigma_v}$.
Then $\sum_{i=1}^m e^{-2\alpha\sigma_{v_i}} = S(\tau_m)$.
Further, by \eqref{eqn:dzt-evol-eqn} the rate of creation of new
vertices at time $t$ is
$|\BP(t)|+B(t)$. Thus, one has
\[
\expec\bigl(d S(t)|\BP(t)\bigr) = e^{-2\alpha t}\bigl(\bigl|\BP
(t)\bigr|+B(t)\bigr)\,dt.
\]
Taking expectations and noting that
$e^{-\alpha t}(|\BP(t)|+ B(t))$ is a martingale gives
\[
\expec\bigl(S(t)\bigr) = 1 + \int_0^t
e^{-\alpha s} \,ds.
\]
{{This completes the proof of part (b), and thus completes the proof of
the lemma. }}
\end{pf}

The next theorem completes the proof of Theorem~\ref
{thm:degree-count-conv}. Before stating the main result, we define some
new constructs which will be used {{in the proof}}. For a bounded
characteristic $\phi$, recall the limit constant $\tilde{m}_\phi
(\infty
)$ in Proposition~\ref{pro:renewal-expec}. In the following theorem, a
key role will be played by the martingale
$ (\tilde{R}_m\dvtx{m\geq0} )$. Recall that this was a
martingale with respect to the filtration $ \{\tilde{\mathcal{F}
}_m\dvtx{m\geq0} \}$.
We shall switch gears and now think about the process
in continuous time. Define $I(t)$ as the set of {{individuals}}
born after time $t$ whose parents were born before time $t$ and note that
%
%
\begin{equation}
\label{eqn:rbpt} \tilde{R}_{|\BP(t)|} {{=}} \sum
_{x\in I(t)} e^{-\alpha\sigma_x}.
\end{equation}
To ease notation, set
%
%
\begin{equation}
\label{eqn:rt-fft-def} R_t:= \tilde{R}_{|\BP(t)|},\qquad \mathcal{F}_t:=
\tilde{\mathcal{F}}_{|\BP(t)|}.
\end{equation}

%
\begin{Theorem}[(Convergence of characteristics)]\label{thm:conv-mart-charac}
For any bounded characteristic that satisfies the recursive
decomposition in \eqref{eqn:charac-recursion}, one has
\[
e^{-\alpha t} Z_\phi(t) \stackrel{a.s.}{\longrightarrow}\tilde
{m}_\phi(\infty)
R_\infty.
\]
Taking $\phi= 1$ and using Proposition~\ref{pro:conts-mp} implies that
{{$R_\infty= W^*$}}, the a.s. limit of the martingale {{$(e^{-\alpha
t}(|\BP(t)|+B(t))\dvtx t\geq0)$}}.
\end{Theorem}

\begin{pf}
First note that Proposition~\ref{prop:repr-mart} implies that $ \{
R_t\dvtx t\geq0 \}$ is an $\mathbb{L}^2$
bounded martingale with respect to the filtration $ \{\mathcal
{F}_t\dvtx
t\geq0 \}$ and
thus $R_t \stackrel{\mathrm{a.s.}}{\longrightarrow}R_\infty$. For a fixed $c
> 0$,
define $I(t,c)$ as the set of vertices born after
time $(t+c)$ whose parents are born before time $t$
and let
%
%
\begin{equation}
\label{eqn:rtc-def} R_{t,c} := \sum_{x\in I(t,c)}
e^{-\alpha\sigma_x}.
\end{equation}
Obviously, $R_{t,c}\leq R_t$. Intuitively, one should expect
$R_{t,c}$ to be small for large $c$. The next lemma
makes this intuition precise. Recall the random
variable $\xi_\alpha= \int_0^\infty e^{-\alpha t}\bolds{\xi}(dt) $ where
{{$\bolds{\xi}= \bolds{\xi}^{v_1}$ denoted the point process
corresponding to
births of offspring of vertex $v_1$.}} {For fixed $c\geq0$, write $\xi
_\alpha(c) := \int_c^\infty e^{-\alpha t} \bolds{\xi}(dt)$. Finally,
define
%
%
\begin{equation}
\label{eqn:xi-alpha-c-def-kc-def} U := \sup_{c\geq
0} e^{c/2}
\xi_\alpha(c),\qquad  A= \mathbb{E}(U),\qquad  K(c) = Ae^\alpha\frac
{e^{-c/2}}{1-\sqrt{e}}.
\end{equation}
The proof below will show that $A< \infty$. Also note that $K(c)\to0$
as $c\to\infty$. }
{{Finally, recall from the proof of Proposition~\ref{pro:conts-mp}
that we defined\break $\lim_{t\to\infty}\exp(-\alpha t) |\BP(t)| = W $.}}
\end{pf}
%
%
\begin{Theorem}\label{prop:rtc-bound}
For any fixed $c>1$, we have
\[
\limsup_{t\to\infty} R_{t,c} \leq K(c) W \qquad \mbox{a.s.},
\]
where $K(c)$ is as in \eqref{eqn:xi-alpha-c-def-kc-def}.
\end{Theorem}

\begin{pf}
The proof uses a variant of the proof used in \cite
{nerman1981convergence}. Let us start by showing that $\mathbb{E}(U) <
\infty$.
First, note that for any fixed $c\geq0$,
\[
e^{c/2}\xi_\alpha(c) \leq\int_c^\infty
e^{t/2} e^{-\alpha t} \bolds{\xi}(dt) \leq\int_0^\infty
e^{t/2} e^{-\alpha t} \bolds{\xi}(dt).
\]
Thus, it is enough to show that $\mathbb{E}(\int_0^\infty e^{t/2}
e^{-\alpha t}
\bolds{\xi}(dt)) < \infty$. By Fubini and integration by parts,
$\mathbb{E}(\int_0^\infty e^{t/2} e^{-\alpha t} \bolds{\xi}(dt)) =
(\alpha-1/2) \int_0^\infty
e^{t/2} e^{-\alpha t} \mu[0,t] \,dt $ where $\mu$ is the intensity
measure of the point process $\bolds{\xi}$. Using Lemma~\ref{lemma:offspring}
shows that for some constant $C< \infty$, we have
\begin{eqnarray*}
\int_0^\infty e^{t/2} e^{-\alpha t}
\mu[0,t] \,dt &\leq& C \int_0^\infty e^{t/2}
e^{-\alpha t} e^{(1-p)t} \,dt\\
& =& C\int_0^\infty
e^{-t/2} \,dt < \infty,
\end{eqnarray*}
by using $\alpha= 2-p$. This completes the proof of finiteness.

Now note that by definition for any $c> 1$
%
%
\begin{eqnarray}
\label{eqn:rtc-bound-one} R_{t,c} &= &\sum_{i=1}^{\lfloor t \rfloor}
\mathop{\sum_{v\dvtx
\sigma_v
\in[i-1,i)}}_{j\dvtx\xi_j^{v}+\sigma_v > t+c } \exp\bigl(-
\alpha\bigl(\xi_j^{v}+\sigma_v\bigr)\bigr) +
\mathop{\sum_{v\dvtx\sigma_v \in[\lfloor
t\rfloor
,t) }}_{ j\dvtx\xi_j^{v}+\sigma_v > t+c } \exp\bigl(-
\alpha\bigl(\xi_j^{v}+\sigma_v\bigr)\bigr)
\nonumber
\\[-8pt]
\\[-8pt]
\nonumber
&\leq&\sum_{i=1}^{\lceil t \rceil}\mathop{ \sum
_{v\dvtx\sigma_v
\in
[i-1,i)}}_{j\dvtx\xi_j^{v}+\sigma_v > t+c } \exp\bigl(-\alpha
\bigl(\xi
_j^{v}+\sigma_v\bigr)\bigr).
\end{eqnarray}
Here, as usual, $\lfloor t \rfloor$ is the largest integer $\leq t$
and $\lceil t \rceil$ is the smallest integer \mbox{$\geq t$}. Analogous to
the definition of $\xi_\alpha(\cdot)$, define for each vertex $v$,
$\xi
_\alpha^v(\cdot)$ using the offspring point process $\bolds{\xi}^v$ of
$v$, namely
\[
\xi_\alpha^v(t):= \int_t^\infty
\exp(-\alpha t) \bolds{\xi}^v(dt)= \sum_{j\dvtx\xi_j^v \geq t}
\exp\bigl(-\alpha\xi_j^v\bigr).
\]
Further analagous to \eqref{eqn:xi-alpha-c-def-kc-def}, for each vertex
$v$ define
\[
U_v(t): = e^{t/2} \xi_\alpha^v(t),\qquad
U_v := \sup_{t\geq0} e^{t/2}
\xi_\alpha^v(t).
\]
Note that
%
%
\begin{equation}
\label{eqn:uv-u-dist} U_v \stackrel{d} {=} U,\qquad U_v(t)
\leq_{\mathrm{st}} U,
\end{equation}
where $U$ is as in \eqref{eqn:xi-alpha-c-def-kc-def} and $\leq_{\mathrm{st}}$
represents stochastic domination. Now for a fixed $i\geq1$ and vertex
$v$ with $\sigma_v \in[i-1, i)$,
\begin{eqnarray*}
\sum_{j\dvtx\xi_j^{v} > t+c-\sigma_v } e^{-\alpha(\xi
_j^{v}+\sigma_v)}& = &e^{-\alpha\sigma_v}
\xi_{\alpha}^v(t+c -\sigma_v) \\
&\leq& e^{-\alpha
(i-1)}
e^{-(t+c - i)/2} U_v(t+c - \sigma_v).
\end{eqnarray*}
Using this in \eqref{eqn:rtc-bound-one} gives
%
%
\begin{equation}
\label{eqn:rtc-bound-two} R_{t,c} \leq\sum_{i=1}^{\lceil t \rceil}
e^{-\alpha(i-1)} e^{-(t+c
- i)/2} \sum_{v\dvtx\sigma_v \in[i-1,i) }
U_v(t+c - \sigma_v).
\end{equation}
To proceed, we will need the following generalization of the strong
law. We paraphrase the following from \cite{nerman1981convergence},
Proposition~4.1.
%
\begin{Proposition}[(Extension of the strong law)]
\label{prop:slln}
Let $ \{n_i\dvtx i\geq1 \}$ be a sequence of integers and let
$(U_{ij}\dvtx1\leq j\leq n_i) $ be a collection of independent random
variables for each fixed $i\geq1$. Suppose that there exists a random
variable $U > 0$ with $\mathbb{E}(U) < \infty$ such that
%
%
\begin{equation}
\label{eqn:xij-bound} |U_{ij}| \leq_{\mathrm{st}} U,\qquad  1\leq j \leq
n_i.
\end{equation}
Further assume
%
%
\begin{equation}
\label{eqn:ni-increase} \liminf_{i\to\infty} \frac
{n_{i+1}}{n_1+\cdots+ n_i} > 0.
\end{equation}
Then
%
%
\begin{equation}
\label{eqn:si-convg} S_i:= \frac{\sum_{j=1}^{n_i} (U_{ij} - \mathbb
{E}(U_{ij}))}{n_i} \stackrel{a.s.}{\longrightarrow}0\qquad \mbox{ as }
i\to
\infty,
\end{equation}
and in fact for any $\varepsilon>0$
%
%
\begin{equation}
\label{eqn:si-io} \sum_{i=1}^\infty
\prob\bigl(|S_i| > \varepsilon\bigr) < \infty.
\end{equation}
\end{Proposition}
Proceeding with the proof, for any interval $\mathcal{I}\subseteq
\mathbb{R}_+$, write
$\BP(\mathcal{I})$ for the collection of vertices born in the
interval $\mathcal{I}$ so
that $\BP(t) \equiv\BP[0,t]$. We will use the above proposition with
$n_i = |\BP[i-1, i)|$ and for each fixed $i$, the collection of random
variables $ \{U_{v}(t+c - \sigma_v)\dvtx v\in\BP[i-1, i)
\}$. This
is a
little subtle since the above proposition is stated for deterministic
sequences but this justified exactly as in the proof of \cite
{nerman1981convergence},
equation
(5.29). First, note that $U_{v}(t+c - \sigma_v)
\leq_{\mathrm{st}} U$ for each fixed $v$. Note that by Proposition~\ref{pro:conts-mp},
\[
\frac{n_{i+1}}{n_1+\cdots+ n_i}:=\frac{|\BP[i,i+1)|}{\BP[0,i)}
\stackrel{\mathrm{a.s.}}{\longrightarrow}e^{\alpha}-1> 0,
\]
as $i\to\infty$, thus \eqref{eqn:ni-increase} is satisfied (almost
surely). Using Proposition~\ref{prop:slln} in \eqref{eqn:rtc-bound-two}
[in particular \eqref{eqn:si-io}] now shows that for any fixed
$\varepsilon> 0$
\[
\limsup_{t\to\infty} R_{t,c} \leq\limsup
_{t\to\infty} \sum_{i=1}^{\lceil t \rceil}
e^{-\alpha(i-1)} e^{-(t+c -i)/2} \bigl(\mathbb{E}(U)+\varepsilon
\bigr) \bigl|\BP[i-1,i)\bigr|.
\]

Using the fact that $e^{-\alpha i} |\BP[i-1,i)|\leq e^{-\alpha i} |\BP
[0,i)| \stackrel{\mathrm{a.s.}}{\longrightarrow}W$, simplifying the above
bound and recalling that we
used $A=\mathbb{E}(U)$, shows that for every fixed $\varepsilon> 0$
\[
\limsup_{t\to\infty} R_{t,c} \leq W (A+\varepsilon)
e^\alpha\frac
{e^{-(c-1)/2}}{1-\sqrt{e}}.
\]
Since $\varepsilon$ was arbitrary, this completes the proof.
\end{pf}
\noqed\end{pf}

\begin{pf*}{Completing the proof of Theorem~\ref
{thm:conv-mart-charac}}
%
Recall that we are dealing with bounded characteristics, that is,
$\Vert\phi\Vert_\infty< C$ for some constant $C$. Without loss of
generality, let $C=1$. We shall show that there exists a constant
$\kappa$ such that for all $\varepsilon> 0$,
%
%
\begin{equation}
\label{eqn:dist-small} \limsup_{t\to\infty} \bigl|e^{-\alpha t}
Z_\phi(t) - \tilde{m}_\phi(\infty) R_\infty\bigr| \leq
\varepsilon(W+ 2\kappa R_\infty).
\end{equation}
Since this is true for any arbitrary $\varepsilon$, this completes the proof.
Fix $\varepsilon> 0$. First, choose $c$ large such that {{the bound in
Theorem~\ref{prop:rtc-bound} satisfies $K(c) <\varepsilon$}}.
Next, for fixed $s> 0$, define {{the truncated characteristic}}
$\phi
_s$ as
%
%
\begin{equation}
\label{eqn:trunc} \phi_s(u) = \cases{ %
\phi(u),&\quad
$u \leq s,$
\vspace*{2pt}\cr
0,&\quad $u> s.$}
\end{equation}
When the branching process is counted by this characteristic, the
contribution of all vertices whose age is more than $s$ is zero.
{{One}} can view {{this}} as a characteristic
used to count ``young'' vertices.
The limit constant for this characteristic by Proposition~\ref
{pro:renewal-expec} is
\[
\tilde{m}_{\phi_s}(\infty) = \int_0^s
e^{-\alpha u} \expec\bigl(\phi(u)\bigr) \,du,
\]
{{where}} $\phi$ is the original characteristic. Note that $\tilde
{m}_{\phi_s}(\infty) \to\tilde{m}_\phi(\infty)$ as the truncation
level $s\to\infty$. Further, writing $\phi^\prime= \phi- \phi_s$,
we can view $\phi^\prime$ as the characteristic counting scores
for ``old'' vertices (vertices of age greater than $s$). With this
notation, we have $Z_\phi(u) = Z_{\phi_s}(u)+ Z_{\phi^\prime}(u)$.

Define
\[
\tilde{m}_{\phi_s}(u) = e^{-\alpha u}\expec\bigl(Z_{\phi_s}(u)
\bigr),\qquad u\geq0.
\]
{{Now choose $s > c$ large enough with $e^{-\alpha s} < \varepsilon$}} such
that for all $u> s-c$
one has $e^{-\alpha s} < \varepsilon$,
$|\tilde{m}_{\phi_s}(\infty) - \tilde{m}_\phi(\infty)| <
\varepsilon$, and
$|\tilde{m}_{\phi_s}(u)-\tilde{m}_{\phi_s}(\infty)|< \varepsilon$.
The constructs $s$ and $c$ shall remain fixed for the rest of the argument.


Let us understand $Z_{\phi_s}(\cdot)$, the
branching process counted according to the truncated
characteristic. We first
observe that since $\phi_s(u) = 0$ when $u> s$, {{this implies that}}
for any time $t>s$, vertices born before time $t-s$ (old vertices)
do not contribute to $Z_{\phi_s}(t)$. {{Define}} $I(t-s)$ {{as the
collection}} of individuals born after time $t-s$ whose parents were
born before time $t$. Then $Z_{\phi_s}(t)$ decomposes as
\[
Z_{\phi_s}(t) = \sum_{v\in I(t-s)}
Z^v_{\phi_s}(t-{{\sigma_v}}),
\]
where $Z^v_{\phi_s}(t-{{\sigma_v}})$ are the contributions to
$Z_{\phi
_s}(t)$ by the descendants of a vertex $v$ born in the interval
$[t-s,t]$. {{Note that by construction, the parent of such a vertex
$v$}} belongs to $\BP(t-s)$.
Further, recall that {{in the definition of $R_{t,c}$ in \eqref
{eqn:rtc-def} we used}} $I(t-s, c)$ for the set of vertices born after
time $(t-s+c)$ whose parents are born before time $t-s$. Then we can
further decompose the above sum as
\[
Z_{\phi_s}(t) = \sum_{x\in I(t-s)\setminus I(t-s,c)}
Z^x_{\phi
_s}(t-\sigma_x)+ \sum
_{x\in I(t-s, c)} Z^x_{\phi_s}(t-\sigma_x).
\]
To simplify notation, write $\mathcal{N}(t-s,c) = I(t-s)\setminus I(t-s,c)$,
that is, the set of individuals born in the interval $[t-s, t-s+c]$ to
parents who were born before time $t-s$. Then we can decompose the
difference as a telescoping sum:
%
%
\begin{equation}
\label{eqn:telescope} e^{-\alpha t} Z_\phi(t) - {{\tilde{m}_\phi}}(
\infty) R_\infty: = \sum_{j=1}^7
E_j(t).
\end{equation}
The definition of these {{seven}} terms $ \{E_i(t)\dvtx1\leq
i\leq7 \}$
are as follows:
\begin{longlist}[(a)]
\item[(a)]$E_1(t)$ is defined by setting
\[
E_1(t)=e^{-\alpha t}Z_{\phi^\prime}(t),\qquad  t\geq0.
\]
Observe that for $E_1(t)$, the only vertices that
contribute are those with age greater than $s$
(since $\phi^\prime(u) = 0$ for $u< s$). In particular, $E_1(t) =
e^{-\alpha t}Z_{\phi^\prime}(t) \leq e^{-\alpha t} |\BP(t-s)| $. Thus,
by Proposition~\ref{pro:conts-mp}, one has $\limsup_{t\to\infty} E_1(t)
\leq \break e^{-\alpha s} W\leq\varepsilon W$ a.s. by choice of $s$.
\item[(b)]$E_2(t)$ is defined by setting
\[
E_2(t):= \sum_{\sss x\in\mathcal{N}(t-s,c)}e^{-\alpha\sigma
_x}
\bigl[e^{-\alpha(t-\sigma_x)}Z^x_{\phi_s}(t-\sigma_x) -
\tilde{m}_{\phi_s}(t-\sigma_x) \bigr].
\]
Note that since in the above sum $x\in\mathcal{N}(t-s,c) $, thus
$\sigma_x >
t-s$. Thus,
\[
\bigl|E_2(t)\bigr| \leq e^{-\alpha(t-s)} \bigl|\mathcal{N}(t-s,c)\bigr| \frac{\sum
_{x\in\mathcal{N}(t-s,c)
} e^{-\alpha(t-\sigma_x)}Z^x_{\phi_s}(t-\sigma_x) - \tilde
{m}_{\phi
_s}(t-\sigma_x)}{|\mathcal{N}(t-s,c)|}.
\]
{{For $E_2(t)$, $\mathcal{N}(t-s,c)$ consists of all children of
parents in
$\BP(t-s)$ that are born in the interval $[t-s, t-s+c]$. Thus,
$|\mathcal{N}
(t-s,c)| \leq\BP(t-s+c)$. In particular, $\limsup_{t\to\infty}
e^{-\alpha(t-s)}|\mathcal{N}(t-s,c)| \leq We^{\alpha c}$.}} Further,
each of the
individuals in $\BP(t-s)$ reproduce at rate at least $1$. One can check
by the strong law of large numbers that
$\liminf_{t\to\infty} |\mathcal{N}(t-s, c)|/|\BP(t-s)|\geq c$
almost surely.
Finally, the terms in the summand (conditional on $\BP(t-s)$) are
independent random variables and each such term in the sum looks like
$X-\expec(X)$, where $X$ is stochastically bounded by the random
variable $Z_{\phi_s}(c)$. A strong law of large numbers argument shows
that $\limsup_{t\to\infty} |E_2(t)| = 0$ a.s.
\item[(c)]$E_3(t)$ is defined as
\[
E_3(t):= \sum_{x\in\mathcal{N}(t-s,c)} e^{-\alpha\sigma_x}
\bigl(\tilde{m}_{\phi_s}(t-\sigma_x) - \tilde{m}_{\phi_s}(
\infty) \bigr).
\]
By the choice of $s$ since $t-\sigma_x \geq s -c$, $|\tilde{m}_{\phi
_s}(t-\sigma_x) - \tilde{m}_{\phi_s}(\infty)|\leq\varepsilon$.
Thus, one has
$|E_3(t)|\leq\varepsilon R_t$. Letting $t\to\infty$, one gets
$\limsup_{t\to
\infty} |E_3(t)| \leq\varepsilon R_\infty$ a.s.
\item[(d)] $E_4(t)$ is defined as
\[
E_4(t):= \tilde{m}_{\phi_s}(\infty) \biggl(\sum
_{x\in\mathcal{N}(t-s,c)} e^{-\alpha\sigma_x}- R_{t-s} \biggr).
\]
For $E_4(t)$, we have $| (\sum_{x\in\mathcal{N}(t-s,c)} e^{-\alpha
\sigma
_x}- R_{t-s} )|= R_{t-s,c}$. Thus,
\[
\limsup_{t\to\infty} E_4(t) \leq\tilde{m}_{\phi_s}(
\infty) K(c) W \leq\tilde{m}_{\phi}(\infty)\varepsilon W,
\]
almost surely by {{Theorem~\ref{prop:rtc-bound} for the asymptotics of
$R_{t,c}$. Here, we have used $\tilde{m}_{\phi_s}(\infty) \leq
\tilde{m}_{\phi}(\infty)$ and that our choice of $c$ guarantees
$K(c) <\varepsilon$}}. To ease notation for the rest of the proof, let
$\kappa$ be a
constant chosen such that $\max(\sup_{u,s\geq0}(\tilde{m}_{\phi
_s}(u)),\tilde{m}_{\phi}(\infty)) < \kappa$. The uniform
boundedness of
$\phi$ guarantees that this can be done. By choice, $\kappa$ is
independent of $s,u$. Thus, the bound for the fourth term simplifies to
$\limsup_{t\to\infty} E_4(t) \leq\kappa\varepsilon W$.
\item[(e)]$E_5(t)$ is defined by setting $E_5(t):= \tilde{m}_{\phi
_s}(\infty
)(R_{t-s}-R_\infty)$. Since\break $R_{t-s}\stackrel{\mathrm{a.s.}}{\longrightarrow
}R_\infty$, $E_5(t)\stackrel{\mathrm{a.s.}}{\longrightarrow}0$.
\item[(f)]$E_6(t)$ is defined by setting $E_6(t) := R_\infty(\tilde
{m}_{\phi_s}(\infty) - \tilde{m}_\phi(\infty))$. By choice of~$s$,
$|E_6(t)|\leq\varepsilon R_\infty$.
\item[(g)]$E_7(t)$ is defined by setting
%
%
\begin{eqnarray}
\label{eqn:e7-def} E_7(t)&:=& e^{-\alpha t} \sum
_{v\in I(t-s,c)} Z_{\phi_s}^v(t-\sigma_v)\nonumber
\\
&=& \sum_{v\in I(t-s,c)} {{e^{-\alpha\sigma_v}}} \bigl(\exp
\bigl(-\alpha(t-\sigma_v)\bigr) Z_{\phi_s}^v(t-
\sigma_v) -\tilde{m}_{\phi
_s}(t-{{\sigma_v}})
\bigr)
\\
&&{} + \sum_{v\in I(t-s,c)} \exp(-\alpha t)
\tilde{m}_{\phi
_s}(t-{{\sigma_v}}).
\nonumber
\end{eqnarray}
Using the strong law of large numbers and arguing as in (b) shows that
the first term goes to zero as $t\to\infty$ a.s. 
Using the constant $\kappa$ defined in (d) above we get
\[
\sum_{v\in I(t-s,c)} \exp(-\alpha t) \tilde{m}_{\phi_s}(t-
\sigma_x) \leq\kappa\sum_{ v\in I(t-s,c)} \exp(-
\alpha\sigma_x)= \kappa R_{t-s,c}.
\]
Using {{Theorem}} \ref{prop:rtc-bound} and the choice of $c$ and
letting $t\to\infty$, we get
\[
\limsup_{t\to\infty} E_7(t) \leq\varepsilon\kappa
R_\infty\qquad \mbox{a.s.}
\]
\end{longlist}

Combining all these bounds, one finally arrives at
\[
\limsup_{t\to\infty} \bigl|e^{-\alpha t} Z_\phi(t) -
\tilde{m}_\phi(\infty) R_\infty\bigr| \leq\varepsilon
\bigl(W+2R_\infty+\kappa(W+R_\infty)\bigr)\qquad \mbox{a.s.}
\]
Since $\varepsilon> 0$ was arbitrary, this completes the proof.
\end{pf*}

\subsection{Time of first birth asymptotics}
\label{sec:height}
For a rooted tree with root $\rho$ (here $\rho=v_1$), there is a
natural notion of a generation of a vertex $v$. This is defined as the
number of edges on the path between $v$ and $\rho$. Thus, $\rho$
belongs to generation zero, all the neighbors of $\rho$ belong to
generation one, and so forth. The aim of this section is to define a
modified notion of generation in $\BP(t)$, {{owing to the fact that
the surgery operation as constructed in Section~\ref{sec:equiv} that
sets up a method to go from the continuous time model to the discrete
time model implies that the object of study are the number of edges to
the closest red vertex on the path to the root $v_1$.}} For each fixed
$k$, we shall define stopping times $\operatorname{Bir}(k)$
representing the first
time an individual in modified generation $k$ is born into the process
$\BP(\cdot)$. We study asymptotics of $\operatorname{Bir}(k)$ as
$k\to\infty$.
In the
next section, we use these asymptotics to understand height asymptotics
for the Superstar model.

Fix $t> 0$. For each vertex $v\in\BP(t)$, let $r(v)$ denote the first
red vertex on the path from $v$ to the original progenitor of the
process $\BP(\cdot)$, namely $v_1$. If $v$ is a red vertex then
$r(v)=v$. Let $d(v)$ be the number of edges on the path between $v$ and
$r(v)$ so that $d(v) = 0$ if $v$ is a red vertex.

Fix $k\geq1$. Let $\operatorname{Bir}(k)$ denote the stopping times
\[
\operatorname{Bir}(k) = \inf\bigl\{t> 0\dvtx\exists v \in\BP(t),
d(v) = k \bigr\}.
\]
In other words, $\operatorname{Bir}(k)$ is the first time that there
exists a red
vertex in $\BP(t)$ such that the subtree consisting of all blue
descendants of this vertex and rooted at this red vertex has an
individual in generation $k$. Here, we use $\operatorname{Bir}$ to
remind the reader
that this is the time of the first \emph{birth} in a particular
generation. The next theorem proves asymptotics for these stopping times.

%
\begin{Theorem}\label{thm:first-birth}
Let ${{\lamb}}(\cdot)$ be the Lambert function \cite
{corless1996lambertw}. We have
\[
\frac{\operatorname{Bir}(k)}{k} \stackrel
{a.s.}{\longrightarrow}\frac{{{\lamb
}}(1/e)}{1-p}\qquad \mbox{as } k\to\infty.
\]
\end{Theorem}

\begin{pf}
Given any rooted tree
$\mathcal{T}
$ and $v\in\mathcal{T}$, we shall let $G(v)$ denote the generation of this
vertex in $\mathcal{T}$. Write $\BP_b^{v_1}(\cdot)$ for the subtree
consisting
of all blue descendants of the original progenitor ${v_1}$ and rooted
at $v_1$. In distribution, this is just a single type continuous time
branching process where each vertex has the same distribution as the
process $\Yu_{1-p}(\cdot) - 1$. Further, let
\[
{\operatorname{Bir}}^*(k) = \inf\bigl\{t\dvtx\exists v\in\BP_b^{v_1}(t),
G(v) = k \bigr\}.
\]
In words, this is the time of first birth of an individual in
generation $k$ for the branching process $\BP_b^{v_1}(\cdot)$. From the
definitions of $\operatorname{Bir}(k), \operatorname{Bir}^*(k)$, we
have $\operatorname{Bir}(k) \leq\operatorname{Bir}^*(k)$.

Much is know about the time of first birth of a single type
supercritical branching process, in particular %
implies that for $\BP_b^{v_1}(\cdot)$, there exists a limit constant
$\beta$ such that
\[
\operatorname{Bir}^*(k)/k \stackrel{\mathrm{a.s.}}{\longrightarrow}\beta.
\]
Here, $\beta$ can be derived as follows. Write $\mu_b$ for the expected
intensity measure of the blue offspring, that is, as in Lemma~\ref
{lemma:offspring}
\[
\mu_b\bigl([0,t]\bigr) = \expec\bigl(c_B[v_1,
t]\bigr) = e^{(1-p)t} -1,\qquad  t\geq0.
\]
For $\theta> 0$, let
\[
\Phi(\theta):= \expec\biggl(\int_0^\infty
e^{-\theta t} c_B(v_1,dt)\biggr), \qquad \theta\in\mathbb{R}.
\]
It is easy to check that this is finite only for $\theta> 1-p$ since
\[
\Phi(\theta)= \theta\int_0^\infty
e^{-\theta t} \mu_b\bigl([0,t]\bigr) \,dt = \frac
{1-p}{\theta-(1-p)}.
\]
For $a> 0$, define
%
%
\begin{equation}
\label{eqn:Lambda-def} \Lambda(a):= \inf\bigl\{\Phi(\theta
)e^{\theta a}\dvtx\theta
\geq1-p \bigr\} = (1-p)a e^{(1-p)a+1}.
\end{equation}
Then by \cite{kingman-age-depen}, Theorem~5, the limit constant $\beta$
is derived as
%
%
\begin{equation}
\label{eqn:beta-def} \beta= \sup\bigl\{a> 0\dvtx\Lambda(a) < 1
\bigr\}.
\end{equation}
From this, it follows that $\beta= {{\lamb}}(1/e)/(1-p)$ where $
{{\lamb}}(\cdot)$ is the Lambert function. Then we have
\[
\limsup_{k\to\infty}\frac{\operatorname{Bir}(k)}{k} \leq\lim
_{k\to\infty
}
\frac{\operatorname{Bir}
^*(k)}{k} \stackrel{\mathrm{a.s.}}{\longrightarrow}\frac{W(1/e)}{1-p}.
\]
This gives an upper bound in Theorem~\ref{thm:first-birth}. Lemma~\ref
{lemma:lb-bl} proves a lower bound and completes the proof.

%
\begin{Lemma}
\label{lemma:lb-bl}
Fix any $\varepsilon>0$ and let $\beta= {{\lamb}}(1/e)/(1-p)$ be the
{{asserted}} limit constant. Then
\[
\sum_{l=1}^\infty\prob\bigl(\operatorname{
Bir}(l) < (1-\varepsilon) \beta l \bigr) < \infty.
\]
Thus, one has $\liminf_{l\to\infty} \operatorname{Bir}(l)/l\geq
\beta$ a.s.
\end{Lemma}
\begin{pf} For ease of notation, for the rest of this proof we shall
write $t_\varepsilon(l) = (1-\varepsilon)\beta l$. In the full
process $\BP(\cdot)$,
two processes occur simultaneously:

(a) New ``roots'' (red vertices) are created. Recall that we used
$R(\cdot)$ for the counting process for the number of red roots.

(b) The blue descendants of each new root have the same distribution
as a single type continuous time branching process with offspring
process have the same distribution as the process $\Yu_{1-p}(\cdot)-1$.

Fix $l\geq2$ and suppose a new red vertex $v$ was created at some time
$\sigma_v < t_\varepsilon(l)$. Let $\BP_b^v(\cdot)$ denote the
subtree of blue
descendants of $v$. Let $\operatorname{Bir}^*(v,l)> \sigma_v$ be the
time of creation
of the first blue vertex in generation $l$ for subtree $\BP_b^v(\cdot
)$. Now $\operatorname{Bir}(l) < t_\varepsilon(l)$ if and only if
there exists a red vertex
$v$ born before $t_\varepsilon(l)$ such that the subtree of blue descendants
of this vertex has a vertex in generation $l$ by this time. For a fixed
red vertex $v \in\BP(\cdot)$, write $A_v(l)$ for this event. Since
$\operatorname{Bir}^*(v,l) -\sigma_v \stackrel{d}{=} \operatorname
{Bir}^*(l)$, conditional on
$\BP
(\sigma_v)$ one has
\[
\prob\bigl(A_v(l)|\BP(\sigma_v) \bigr) = \prob\bigl(
\operatorname{Bir}^*(l) \leq t_\varepsilon(l) - \sigma_v\bigr).
\]

Fix $0< s< (1-\varepsilon)\beta l$. Then for $\theta> 1-p$, Markov's
inequality implies
\[
\prob\bigl(\operatorname{Bir}^*(l) < (1-\varepsilon)\beta l - s
\bigr) \leq e^{\theta((1-\varepsilon
)\beta
l -s)} \expec
\bigl[e^{-\theta\operatorname{Bir}^*(l)}\bigr].
\]
One of the main bounds of Kingman (\cite{kingman-age-depen},
equation (2.5), Theorem~1) is $\expec[e^{-\theta\operatorname
{Bir}^*(l)} ] \leq(\Phi
(\theta))^l$. Thus, we get
%
%
\begin{equation}
\label{eqn:prob-upb} \prob\bigl(\operatorname{Bir}^*(l) <
(1-\varepsilon)\beta l - s \bigr) \leq\bigl[\Phi(
\theta)e^{\theta
(1-\varepsilon)\beta}\bigr]^l e^{-\theta s}.
\end{equation}
By the definition of $\beta$,
\[
\Lambda_\varepsilon: =\Lambda\bigl(\beta(1-\varepsilon)\bigr) :=
\inf\bigl\{\Phi(
\theta)e^{\theta(1-\varepsilon)\beta}\dvtx\theta> 1-p \bigr\} < 1,
\]
where $\Lambda$ is as in \eqref{eqn:Lambda-def}. It is easy to check
that the minimizer occurs at
\[
\theta_\varepsilon= 1-p+ \frac{1}{(1-\varepsilon)\beta}.
\]
The final probability bound we shall use is
%
%
\begin{equation}
\label{eqn:pb-bd-ldp} \prob\bigl(\operatorname{Bir}^*(l) <
(1-\varepsilon)\beta l - s \bigr) \leq[
\Lambda_\varepsilon]^l e^{-\theta
_\varepsilon s}.
\end{equation}
Let $N_l^\varepsilon$ be the number of red vertices born before time
$t_l(\varepsilon)$
whose trees of blue descendants $\BP_b^v(\cdot)$ have
at least one vertex in generation $l$ by time $t_\varepsilon(l)$.
Obviously, $\prob(\operatorname{Bir}(l) < (1-\varepsilon) \beta l)
\leq\expec(N_l^\varepsilon)$.
Conditioning on the times of birth of red vertices, one gets
\begin{eqnarray*}
\expec\bigl(N_l^\varepsilon\bigr) &\leq&\int
_0^{t_\varepsilon(l)} [\Lambda_\varepsilon]^l d
\expec\bigl(R(s)\bigr) \qquad\mbox{using equation }\eqref{eqn:pb-bd-ldp},
\\
&=&p[\Lambda_\varepsilon]^l\int_0^{t_\varepsilon(l)}
e^{-(\theta_\varepsilon-q)s} \,ds \qquad\mbox{using Lemma~\ref{lemma:offspring}.}
\end{eqnarray*}
Simplifying, we get for all $l\geq2$, $\expec(N_l^\varepsilon) \leq C
[\Lambda
_\varepsilon]^l$ for a constant $C$. Thus,
\[
\sum_{l=1}^\infty P\bigl(\operatorname{Bir}(l) < (1-\varepsilon)
\beta l \bigr) < \infty.
\]
\upqed\end{pf}

\section{Proofs of the main results}
\label{sec:comp-proof}
Recall the equivalence created by the surgery operation between the
Superstar model and the two-type branching process as established in
Section~\ref{sec:equiv}. We shall use this equivalence and the proven
results on $\BP(\cdot)$ in Section~\ref{sec:proofs} to complete the
proof of the main results. We record the following fact about the
asymptotics for the stopping times $\tau_n$.
%
%
\begin{Lemma}[(Stopping time asymptotics)]\label{lemma:stop-times}
The stopping times $\tau_n$ satisfy
\[
\tau_n - \frac{1}{2-p} \log{n} \stackrel{a.s.}{\longrightarrow
}-\frac{1}{2-p}
\log{W}.
\]
\end{Lemma}
\begin{pf} Proposition~\ref{pro:conts-mp} proves that $|\BP
(t)|e^{-(2-p)t} \stackrel{\mathrm{a.s.}}{\longrightarrow}W$. Thus\break 
$ne^{-(2-p)\tau_n}\stackrel{\mathrm{a.s.}}{\longrightarrow}W$.
\end{pf}

Let us now start by proving the main results.
We note that Theorem~\ref{SSLLN} is obvious since the degree of the
superstar is given by $R(\tau_n) = \sum_{i=1}^n\ind\{v_i \mbox
{ is red}\}$,
the total number of red vertices and $(\ind\{v_i \} \mbox
{ is red})_{i\geq1}$ is an i.i.d. sequence with Bernoulli $p$ as the
marginal distribution. We now prove the remaining results using the
correspondence between the continuous time and discrete time processes.

\subsection{Proof of the degree distribution strong law}
In this section, we shall prove Theorem~\ref{thm:DegreeDist}. Since
$G_{n+1}$ is a connected tree, every vertex has degree at least one.
Recall that $c_B(v,t)$ denotes the number of blue children of vertex
$v$ by time~$t$. Write $\deg(v, G_{n+1})$ for the degree of a vertex in
$G_{n+1}$. The surgery operation implies that for any nonsuperstar vertex
%
%
\begin{equation}
\label{eqn:degree-blue-rel} \deg(v, G_{n+1}) = c_B(v,
\tau_n) + 1.
\end{equation}
Fixing $k\geq0$, the number of nonsuperstar vertices with degree
exactly $k+1$ is the same as the number of vertices in $\BP(\tau_n)$
that have exactly $k$ blue children. Recall that we used $Z_{\geq
k}(t)$ for the number of vertices in $\BP(t)$ that have at least $k$
blue children. Proposition~\ref{pro:conts-mp}, showed that the total
number of vertices $|\BP(t)|$ satisfies
%
%
\begin{equation}
\label{eqn:bp-convg-res} e^{-(2-p)t}\bigl|\BP(t)\bigr| \stackrel
{\mathrm{a.s.}}{\longrightarrow}\frac
{W^*}{(2-p)} \qquad\mbox{as } t
\to\infty.
\end{equation}
Theorem~\ref{thm:degree-count-conv} showed that
\[
e^{-(2-p)t}Z_{\geq k}(t) \stackrel{\mathrm{a.s.}}{\longrightarrow}k!\prod
_{i=1}^k \biggl(i+\frac
{2-p}{1-p}
\biggr)^{-1} \frac{W^*}{2-p}.
\]
Thus, writing $p_{\geq k}(t) = Z_{\geq k}(t)/\BP(t)$ for the proportion
of vertices with degree $k$, Theorem~\ref{thm:degree-count-conv}
implies one has
\[
p_{\geq k}(t) \stackrel{\mathrm{a.s.}}{\longrightarrow}k!\prod_{i=1}^k
\biggl(i+\frac{2-p}{1-p} \biggr)^{-1}:= p_{\geq k}(\infty)\qquad
\mbox{as $t\to\infty$.}
\]
Now let $k\geq1$. Writing $N_{\geq k}(n)$ for the number of vertices
with degree at least $k$ in $G_{n+1}$, one has $N_{\geq
k}(n)/n\stackrel{\mathrm{a.s.}}{\longrightarrow}
p_{\geq k-1}(\infty)$ as $n\to\infty$. Thus, the proportion of vertices
with degree exactly $k$ converges to $p_{\geq k-1}(\infty) - p_{\geq
k}(\infty) = \nu_{\mathrm{SM}}(k)$. This completes the proof.

\subsection{Proof of maximal degree asymptotics}
\label{sec:max-deg-asympt}
The aim of this is to prove Theorem~\ref{thm:DegreeEvolutionLaw}. We wish
to analyze the maximal nonsuperstar degree that we wrote as
\[
\Upsilon_n = \max\bigl\{\deg(v_i, G_{n+1})
\dvtx1\leq i\leq n \bigr\}.
\]
The plan will be as follows: we will first prove the simpler assertion
of convergence of the degree of vertex $v_k$ for fixed $k\geq1$. Then
we shall show that given any $\varepsilon> 0$, we can choose $K$ such
that for
large $n$, the maximal degree vertex has to be one of the first $K$
vertices $v_1, v_2,\ldots, v_K$ with probability greater than
$1-\varepsilon$.
This completes the proof.

Fix $k\geq1$. Recall from \eqref{eqn:degree-blue-rel} that
$\deg(v_k, G_{n+1}) = c_B(v_k, \tau_n)+1$ where $c_B(v_k, t)$
are the number of blue vertices born to vertex $k$ by time $t$.
Recall that $c_B(v_k, t)$ is a Yule process of rate $1-p$ started
at time $\tau_k$ (i.e., at the birth of vertex $v_k$). By Lemma~\ref
{lem:yule-proc},
%
%
\begin{equation}
\label{eq:bpc} \frac{ c_B(v_k, t)}{e^{(1-p) (t-\tau_k)}} \stackrel
{\mathrm{a.s.}} {\longrightarrow}
W^\prime_k,
\end{equation}
where $W^\prime_k$ is an exponential random variable with mean one.
Write $\gamma= (1-p)/(2-p)$ and let $\Delta_k = e^{-(1-p)\tau_k}
W'W^{-\gamma}$. Using \eqref{eqn:bp-convg-res} and \eqref{eq:bpc},
we have
\begin{eqnarray*}
n^{-\gamma}\deg(v_k,G_{n+1}) &=& \frac{ c_B(v_k, \tau_{n-1})+1}{e^{(1-p)
(\tau_{n-1}-\tau_k)}}
\biggl(\frac{e^{(2-p)\tau_{n-1}}}{|\BP(\tau_{n-1})|+1} \biggr
)^\gamma e^{-(1-p)\tau_k}
\\
&\stackrel{\mathrm{a.s.}}{\longrightarrow}&W_k'W^{-\gamma}e^{-(1-p)\tau_k}:=
\Delta_k.
\end{eqnarray*}

Now let us prove distributional convergence of the properly normalized
maximal nonsuperstar
degree $\Upsilon_n$. Fix $L> 0$ and let
%
%
\begin{equation}
\label{eqn:} \tilde{M}_n[0,L]:= \max\bigl\{\deg(v_k,
G_{n+1})\dvtx\tau_k \leq L \bigr\}.
\end{equation}
In other words, this is the largest degree in $G_{n+1}$ amongst all
vertices born before time $L$ in $\BP(\cdot)$. The convergence of the
degree of $v_k$ for any $k\geq1$ implies the next result.
%
%
\begin{Lemma}[(Convergence near the root)]\label{lem:convg-near-root}
Fix any $L> 0$. Then there exists a random variable $\Delta^*[0,L]> 0$
such that
\[
{n^{-\gamma}\tilde{M}_n[0,L]}\stackrel{a.s.} {\longrightarrow}
\Delta^*[0,L],
\]
where $\gamma=(1-p)/(2-p)$.
\end{Lemma}
Now if we can show that with high probability, $\Upsilon_n = \tilde
{M}_n[0,L]$ for large finite $L$ as $n\to\infty$, then we are done.
This is accomplished via the next lemma. Recall that by asymptotics for
the stopping times $\tau_n$ in Lemma~\ref{lemma:stop-times}, given any
$\varepsilon> 0$, we can choose $K_\varepsilon> 0$ such that
%
%
\begin{equation}
\label{eqn:tn-b-bound} \limsup_{n\to\infty}\prob\biggl(\biggl
\llvert
\tau_n - \frac{1}{2-p}\log{n}\biggr\rrvert> K_\varepsilon
\biggr) \leq\varepsilon.
\end{equation}
For any $0< L< t$, let $\BP(L,t]$ denote the set of vertices born in
the interval $(L,t]$. Recall that we used ${v_1}$ for the original
progenitor. For any time $t$ and $v\in\BP(t)$, let $\operatorname
{deg}_v(t) =
c_B(v,t)+1$ denote the degree of vertex $v$ in the Superstar model
$G_{|\BP(t)|+1}$ obtained through the surgery procedure. For fixed $K$
and $L$, let $A_n(K,L)$ denote the event that for some time $t\in
[(2-p)^{-1}\log{n}\pm K]$, there exists a vertex $v$ in $\BP(L,t]$ with
$\deg_v(t) > \deg_{v_1}(t)$.
%
%
\begin{Lemma}[(Maxima occurs near the root)]\label{lem:max-occurs-root}
Given any $K$ and $\varepsilon$, one can choose $L > 0$ such that
\[
\limsup_{n\to\infty}\prob\bigl(A_n(K,L)\bigr) \leq\varepsilon.
\]
In particular, given any $\varepsilon> 0$, we can choose $L$ such that
\[
\limsup_{n\to\infty}\prob\bigl(\Upsilon_n \neq
\tilde{M}_n\bigl([0,L]\bigr)\bigr) \leq\varepsilon.
\]
\end{Lemma}
Deferring the proof of this result note that Lemma~\ref
{lem:convg-near-root} now coupled with the above lemma now shows that
there exists a random variable $\Delta^*$ such that
$\Upsilon_n/n^{\gamma} \convp\Delta^*$.
This completes the proof of Theorem~\ref{thm:DegreeEvolutionLaw}.

\begin{pf*}{Proof of Lemma~\ref{lem:max-occurs-root}} For ease of
notation, write
\[
t_n^- = (2-p)^{-1}\log{n} -K, \qquad t_n^+ =
(2-p)^{-1}\log{n} +K.
\]
Since the degree of any vertex is an increasing process it is enough to
show that we can choose $L = L(K,\varepsilon)$ such that as $n\to
\infty$, the
probability that there is some vertex born in the time interval $[L,
t_n^+]$ whose degree at time $t_n^+$ is larger than the degree of the
root ${v_1}$ at time $t_n^-$ is smaller than $\varepsilon$. Let
$M_{[L,t_n^+]}(t_n^+)$ denote the maximal degree by time $t_n^+$ of all
vertices born in the interval $[L,t_n^+]$. Then for any constant $C> 0$
\begin{eqnarray*}
\prob\bigl(A_n(K,L)\bigr) &\leq&\prob\bigl( \bigl\{
\deg_{v_1}\bigl(t_n^-\bigr) < C n^\gamma\bigr\}
\cup\bigl\{M_{[L,t_n^+]}\bigl(t_n^+\bigr) > C n^{\gamma
}
\bigr\} \bigr)
\\
&\leq&\prob\bigl({\deg_{v_1}\bigl(t_n^-\bigr) < C
n^\gamma} \bigr)+ \prob\bigl({ M_{[L,t_n^+]}\bigl(t_n^+
\bigr) > C n^\gamma} \bigr).
\end{eqnarray*}
Since the offspring process of ${v_1}$ has the same distribution as a
rate $(1-p)$ Yule process
\[
e^{-(1-p)t_n^-}\deg_{v_1}\bigl(t_n^-\bigr) =
e^{(1-p)K/2}\frac{\deg
_{v_1}(t_n^-)}{n^\gamma} \stackrel{\mathrm{a.s.}}{\longrightarrow}W_{v_1},
\]
where $W_{v_1}$ has an exponential distribution with mean one. Thus,
for a fixed $K$, we can choose $C = C(\varepsilon)$ large enough such that
\[
\limsup_{n\to\infty} \prob\bigl({\deg_{v_1}
\bigl(t_n^-\bigr) < C n^{\gamma
}} \bigr) \leq\varepsilon/2.
\]
Thus, for a fixed $\varepsilon, C, K$, it is enough to choose $L$
large such that
\[
\limsup_{n\to\infty} \prob\bigl({ M_{[L,t_n^+]}
\bigl(t_n^+\bigr) > C n^{\gamma
}} \bigr) \leq\varepsilon/2.
\]
Without loss of generality, we shall assume throughout that
$L_\varepsilon$
and $ t_n^+$ are integers. For any integer $L_\varepsilon< m < t_n^+-1
$, let
$M_{[m,m+1]}(t_n^+)$ denote the maximum degree by time $t_n^+$ of all
vertices born in the interval $[m,m+1]$. Then
\[
M_{[L,t_n^+]}\bigl(t_n^+\bigr) = \max_{L\leq m \leq t_n^+-1}
M_{[m,m+1]}\bigl(t_n^+\bigr).
\]
Let $|\BP[m,m+1]|$ denote the number of vertices born in the time
interval $[m,m+1]$. Since for a vertex born at some time $s< t_n^+$,
the degree of the vertex at time $t_n^+$ has distribution
$\operatorname
{Yu}_{1-p}(t_n^+ -s)$, an application of the union bound yields
\[
\prob\bigl({ M_{[L,t_n^+]}\bigl(t_n^+\bigr) > C
n^{\gamma}} \bigr) \leq\sum_{m=L}^{t_n^+-1}
\expec\bigl(\bigl|\BP[m,m+1]\bigr|\bigr) \prob\bigl(\operatorname{Yu}_{1-p}
\bigl(t_n^+ -m\bigr) > C n^{\gamma}\bigr).
\]
Now $\expec(\BP[m,m+1])\leq\expec(|\BP(m+1)|)$. By
Proposition~\ref
{pro:conts-mp}, $\expec(|\BP(t)|) \leq e^{(2-p)t}$. Further by
Lemma~\ref{lem:yule-proc}, for fixed time $s$, a rate $1-p$ Yule
process has
a geometric distribution with parameter $e^{-(1-p)s}$. Thus, we have
%
%
\begin{eqnarray}
\label{eqn:max-degree-different-levels} \prob\bigl({
M_{[L,t_n^+]}\bigl(t_n^+\bigr) > C
n^{\gamma}} \bigr) &\leq&\sum_{m=L}^{t_n^+-1}
Ae^{(2-p)m} \bigl[1-e^{-(1-p)(t_n^+ - m)} \bigr]^{ C
n^\gamma}
\nonumber
\\[-8pt]
\\[-8pt]
\nonumber
&\leq&\sum_{m=L}^{t_n^+-1} Ae^{ ((2-p)m - Ce^{(1-p)(m-K)} )},
\end{eqnarray}
where last inequality follows from the fact that for $0\leq x\leq1$,
$1-x\leq e^{-x}$ and
\[
e^{t_n^+/2} = n^{\gamma}e^{(1-p)K}.
\]
Now choosing $L$ large, one can make the right-hand side of the last
inequality as small as one desires and this completes the proof.
\end{pf*}

\subsection{Proof of logarithmic height scaling}
The aim of this section is to complete the proof of Theorem~\ref
{thm:height}. Let us first understand the relationship between the
distances in $\BP(\tau_n)$ and $G_{n+1}$ due to the surgery operation.
The distance of all the red vertices in $\BP(\tau_n)$ from the
superstar $v_0$ is one. For each blue vertex $v\in\BP(\tau_n)$, let
$r(v)$ denote the first red vertex on the path from $v$ to the root
$v_1$ in $\BP(\tau_n)$. Recall from Section~\ref{sec:height} that
$d(v)$ denoted the number of edges on the path between $v$ and $r(v)$
with $d(v) = 0$ if $v$ was a red vertex. Then the distance of this
vertex from the superstar $v_0$ in $G_{n+1}$ is just $d(v)+1$ since the
vertex needs $d(v)$ steps to get to $r(v)$ that is then directly
connected to $v_0$ in $G_{n+1}$ by an edge. Let $D(u,v)$ denote the
graph distance between vertices $u$ and $v$ in $G_{n+1}$. Since by
convention $d(v)=0$ for all the red vertices, this argument shows that
for all $v\neq v_0\in G_{n+1}$, $ D(v, v_0) = d(v)+1$. In particular,
the height of $G_{n+1}$ is given by
%
%
\begin{equation}
\label{eqn:hhgn-dvn-equiv} \HH(G_{n+1}) = \max\bigl\{d(v)+1\dvtx
v\in\BP(
\tau_n) \bigr\}.
\end{equation}
Now by the definition of $\HH(G_{n+1})$, there is a vertex in $\BP
(\tau
_n)$ such that $d(v)=\HH(G_{n+1})-1$ but no vertex with $d(v) = \HH
(G_{n+1})$. Recall the stopping times $\operatorname{Bir}(k)$, defined
as the first
time a vertex with $d(v)=k$ is born in $\BP(\cdot)$. Thus, we have
%
%
\begin{equation}
\label{eqn:ineq} \operatorname{Bir}\bigl(\HH(G_{n+1})-1\bigr) \leq
\tau_n
\leq\operatorname{Bir}\bigl(\HH(G_{n+1})\bigr).
\end{equation}
Now recall that Theorem~\ref{thm:first-birth} showed that the stopping
times $\operatorname{Bir}(k)$
satisfy
\[
\operatorname{Bir}(k)/k\stackrel{\mathrm{a.s.}}{\longrightarrow}{{\lamb
}}(1/e)/(1-p) \qquad\mbox{as } k\to
\infty.
\]
Dividing \eqref{eqn:ineq} throughout by $\HH(G_{n+1})$ by
Theorem~\ref
{thm:first-birth}
\[
\frac{\operatorname{Bir}({\HH(G_{n+1})-1})}{\HH(G_{n+1})} \stackrel
{\mathrm{a.s.}}{\longrightarrow}
\frac{{{\lamb
}}(1/e)}{1-p},
\]
while by Lemma~\ref{lemma:stop-times} we get
\[
\frac{\tau_n}{\log{n}}\stackrel{\mathrm{a.s.}}{\longrightarrow}\frac{1}{2-p}.
\]
Rearranging shows that
\[
\frac{\HH(G_{n+1})}{\log{n}} \stackrel{\mathrm{a.s.}}{\longrightarrow}\frac
{(1-p)}{{{\lamb}}(1/e)(2-p)}.
\]
This completes the proof.
\end{pf}

\subsection{Extension to the variants of the Superstar model}
\label{sec:extension-proofs}
{{We now describe how the above methodology easily extends to the two
variants described in Section~\ref{sec:discussion}, namely the
superstar linear preferential attachment and the uniform attachment
model (Theorems \ref{thm:superstar-linear} and \ref
{thm:superstar-uniform}). Since the proofs are identical to the
original model, modulo the driving continuous time branching process,
we will not give full proofs but rather describe the continuous time
versions that need to be analyzed to understand the corresponding
discrete model. The surgery operation and the subsequent analysis of
the continuous time model are identical to the original Superstar
model.

For fixed $a> -1$ and $p\in(0,1)$, we write $ \{
G_n^{\operatorname{lin}
}(a,p)\dvtx n\geq1 \}$ for the corresponding family of growing
random trees obtained via following the dynamics of the linear
attachment scheme (see Section~\ref{sec:discussion}). We let $ \{
G_n^{\sss\operatorname{uni}}(p)\dvtx n\geq1 \}$ be the family of random
trees obtained via uniform attachment. Now recall that the analysis of
the superstar preferential attachment model start with the formulation
of a continuous time two type branching process (consisting of red and
blue vertices). One then performs surgery on this two type branching
process at appropriate stopping times $\tau_n$ as defined in \eqref
{eqn:stop-time-def} to obtain the Superstar model. For the two
variants, let us now describe the corresponding continuous time
versions.
\begin{longlist}[(a)]
\item[(a)]\textit{Superstar linear preferential attachment}: We write
$ \{\BP_{\operatorname{lin}}(t) \}_{t\geq0}$ for this branching
process. Here one starts with a single red vertex $v_1$ at time $t=0$.
Each individual lives forever. For any fixed $t\geq0$, each individual
$v\in\BP_{\operatorname{lin}}(t)$ in the branching process
reproduces at rate
\[
\lambda(v,t):= c_B(v,t)+1+a,
\]
where as before $c_B(v,t)$ denotes the number of blue children of
vertex $v$ at time~$t$. Each new offspring is colored red with
probability $p$ and blue with probability $q:=1-p$.
\item[(b)]\textit{Uniform attachment}: Start with a single red vertex
$v_1$ at time $t=0$. Each individual reproduces at rate one and lives
forever. Each new offspring is colored red with probability $p$ and
blue with probability $q:=1-p$. Write $ \{\BP_{\operatorname{uni}}(t)
\}
_{t\geq0}$ for this branching process.
\end{longlist}
Fix $n\geq1$ and recall the stopping time $\tau_n$ from \eqref
{eqn:stop-time-def}, namely the time for the branching process to reach
size $n$. From Section~\ref{sec:equiv}, recall the surgery operation
that takes $\BP(\tau_n)$ to a random tree $\cS_n$ on $n+1$ vertices.
The following proposition which is the general analog of
Proposition~\ref{prop:equiv} showing the equivalence of the continuous
time models
and the discrete time versions. The result is stated for the linear
preferential attachment model, the same result is true using the
corresponding branching process for the uniform attachment model.
%
%
\begin{Proposition} Fix $a> -1$ and $p\in(0,1)$. Let $ \{\BP
_{\operatorname{lin}}(t)\dvtx t\geq0 \}$ be the continuous time two type
branching process constructed as above for the superstar linear
preferential attachment model with parameters $a, p$. The sequence of
trees $ \{\cS_n\dvtx n\geq1 \}$ obtained by performing the
surgery operation on $ \{\BP_{\operatorname{lin}}(\tau_n)\dvtx n\geq1
\}$ has the same distribution as $ \{G_{n+1}(a,p)\dvtx n\geq
1 \}$.
\end{Proposition}
Now recall that in the proof of the original Superstar model, a major
role was played by Proposition~\ref{pro:conts-mp} which showed that
the associated continuous time branching process grew at rate $\exp
((2-p)t)$. This allowed us to make rigorous the following two ideas
(see, e.g., the proof of Corollary~\ref{cor:pk-expression}):
\begin{longlist}[(a)]
\item[(a)] As $t\to\infty$, the age of an individual chosen uniformly at
random from the population has an exponential distribution with rate
$(2-p)$.
\item[(b)] For vertex $v$, let $c_B(v,\sigma_v+t)$ be the number of blue
children $t$ units after being born and note that $ \{c_B(\sigma
_v + t)\dvtx t\geq0 \}$ has the same distribution for any
vertex. Since the number of blue children of a vertex represents the
out-degree in the Superstar model after the surgery operation, using
(i), the limiting degree distribution should be the same as
$1+c_B(v_1,T)$, where $T\sim\exp((2-p))$ independent of $ \{
c_B(v_1,t)\dvtx t\geq0 \}$. Here, we use $v_1$ for convenience
since $\sigma_{v_1} = 0$.
\end{longlist}
The corresponding version of Proposition of \ref{pro:conts-mp} is the
following.
%
%
\begin{Proposition} \label{pro:conts-mp-variants}
%
%
\textup{(a)} Fix $a> -1$ and $p\in(0,1)$. Then there exists a random variable
$W(a,p) > 0$ a.s. such that as $t\to\infty$,
\[
\exp\bigl(-(2-p+a)\bigr)\bigl|\BP_{\operatorname{lin}}(t)\bigr| \stackrel
{a.s.}{\longrightarrow}W(a,p).
\]
\textup{(b)} For the uniform attachment model, for any $p\in(0,1)$ as $t\to
\infty$,
\[
\exp(-t) \bigl|\BP_{\operatorname{uni}}(t)\bigr| \stackrel
{a.s.}{\longrightarrow}W,
\]
where $W\sim\exp(1)$.
%
\end{Proposition}
\begin{pf} We start with part (b).
For the uniform attachment
model, since every individual lives forever and reproduces at rate one,
the process $ \{|\BP_{\operatorname{uni}}(t)|\dvtx t\geq0 \}$ has the
same distribution as a rate one Yule process (see Definition~\ref
{def:yule}). Then the result follows from Lemma~\ref{lem:yule-proc}.

To prove (a), define the process
\[
M(t):= \exp\bigl(-(2-p+a)\bigr) \bigl(\bigl|\BP_{\operatorname
{lin}}(t)\bigr|+ B(t)\bigr),\qquad  t
\geq0,\vadjust{\goodbreak}
\]
where as before $B(t)$ denotes the number of blue individuals in the
population by time $t$. Arguing exactly as in the proof of
Proposition~\ref{pro:conts-mp}, it is easy to check that this process
is a
martingale. The rest of the proof now follows along the same lines as
the proof of Proposition~\ref{pro:conts-mp}.
\end{pf}

The proof of
Theorems \ref{thm:superstar-linear} and \ref{thm:superstar-uniform}
now proceed as in the analysis of the original model. For example, to
show the convergence of the degree distribution for the uniform
attachment model Theorem~\ref{thm:superstar-uniform}, first note that
for any vertex $v$, since this vertex reproduces at rate one and each
new offspring is colored red with probability $p$ and blue with
probability $q=1-p$. Thus, the process counting the number of blue
children $ \{c_B(v_1,t)\dvtx t\geq0 \}$ is a rate $q$
Poisson process. Fix $k\geq1$ and write $Z_{\geq k}(t)$ for the number
of vertices in $\BP_{\operatorname{lin}}(t)$ which have $k$ or more
blue offspring
by time $t$. The analogous version of Theorem~\ref
{thm:degree-count-conv} for the uniform attachment model implies that
\[
\exp(-t) Z_{\geq k}(t) \stackrel{\mathrm{a.s.}}{\longrightarrow}p_{\geq
k}(\infty) W,
\]
where
\[
p_{\geq k}(\infty) = \prob\bigl(c_B(v,T) \geq k\bigr),
\]
where $T\sim\exp(1)$ independent of $c_B(\cdot)$. Now note that
\[
\prob\bigl(c_B(v,T) \geq k\bigr) = \prob\Biggl(\sum
_{i=1}^k \xi_i \leq T\Biggr),
\]
where $ \{\xi_i \}_{i\geq1}$ is a sequence of independent
rate $q$ exponential random variables. Arguing as in the proof of
Corollary~\ref{cor:pk-expression}, we get
\[
\prob\Biggl(\sum_{i=1}^k
\xi_i \leq T \Biggr) = \bigl(\expec\bigl(\exp(-\xi_i)
\bigr)\bigr)^k = \biggl(\frac{q}{q+1} \biggr)^k.
\]
For the maximal degree, note that by Proposition~\ref
{pro:conts-mp-variants} implies that the stopping time $\tau_n$ as in
\eqref{eqn:stop-time-def} for the time the continuous time branching
process grows to be of size $n$ satisfies
\[
\tau_n = \log{n} + O_P(1).
\]
Since for each vertex, its true degree is the number of blue offspring,
as an easy lower bound, the root $v_1$ by time $\tau_n$ should have
degree $\sim(1-p) \log{n}$ (since the process describing the blue
offspring of the root is just a rate $q$ Poisson process). To get that
$\log{n}$ is the correct order for the maximal degree and in
particular the weak law, one argues as in Section~\ref
{sec:max-deg-asympt} [in particular see \eqref
{eqn:max-degree-different-levels}], teasing apart the contribution to
this maximal degree of vertices born at various times. The proof of
Theorem~\ref{thm:superstar-linear} is similar. We omit the details. }}

\begin{appendix}\label{app}
\section*{Appendix}
Below we describe each of the thirteen events and show the
corresponding event specific term.

\begin{itemize}
\item$E=1$: Brazil vs. Netherlands soccer match from the 2010 World
Cup. The term
is ``Brazil'' or ``Netherlands.''

\item$E=2$: Basketball player Lebron James announcement of signing
with the Miami Heat. The term
is ``Lebron.''

\item$E=3$: The 2010 World Cup Kick-Off Celebration Concert. The term
is ``World Cup.''

\item$E=4$: Brazil vs. Portugal soccer match from the 2010 World Cup.
The term is ``Brazil'' or ``Portugal.''

\item$E=5$: Italy vs Slovakia soccer match from the 2010 World Cup.
The term is ``Italy'' or ``Slovakia.''

\item$E=6$: The 2010 BET Awards show. The term is ``BET Awards.''

\item$E=7$: The firing of General Stanly McChrystal by US President
Barack Obama. The term is ``McChrystal.''

\item$E=8$: The 2010 World Cup Opening Ceremony. The term is ``World Cup.''

\item$E=9$: Mexico vs. South Africa soccer match from the 2010 World
Cup. The term is ``Mexico.''

\item$E=10$: England vs. Slovakia soccer match from the 2010 World
Cup. The term is ``England.''

\item$E=11$: Portugal vs. North Korea soccer match from the 2010 World
Cup. The term is ``Portugal.''

\item$E=12$: Roger Federer's tennis match in the first round of the
2010 Wimbledon tournament. The term is ``Federer.''

\item$E=13$: The UN imposing sanctions on Iran. The term is ``Iran.''
\end{itemize}
\end{appendix}

\section*{Acknowledgements}
S. Bhamidi would like to thank the hospitality of the
Statistics department at Wharton where this work commenced.
We thank two referees for detailed comments which lead to a significant
improvement in the readability of the paper.

%
%





\printaddresses
\end{document}